\documentclass[a4paper,11pt,oneside]{article}

\usepackage[top=2cm, bottom=2cm, left=2cm, right=2cm]{geometry}
\usepackage{amsmath,amsthm,amsfonts,amssymb,amscd}
\usepackage{algorithm}
\usepackage{algorithmic}
\usepackage{hyperref}
\usepackage{authblk}
\theoremstyle{plain} % default
\newtheorem{theorem}{Theorem}[section]

\newtheorem{corollary}[theorem]{Corollary}
\newtheorem{remark}[theorem]{Remark}
\newtheorem{definition}[theorem]{Definition}
\newtheorem{example}[theorem]{Example}

\theoremstyle{definition} %

\theoremstyle{remark} %

\usepackage{amssymb}
\usepackage{xcolor}
\usepackage{xargs}
\usepackage[colorinlistoftodos,prependcaption,textsize=tiny]{todonotes}
\usepackage{comment}

\begin{document}

\newcommand{\R}{{\mathbf{R}}}
\newcommand{\V}{{\mathbf{V}}}
\newcommand{\W}{{\mathbf{W}}}
\newcommand{\Z}{{\mathbf{Z}}}
\newcommand{\Dsf}{\mathsf{D}}
\newcommand{\magenta}{\color{magenta}}
\newcommand{\black}{\color{black}}
\newcommand{\blue}{\color{blue}}
\newcommand{\red}{\color{red}}
\newcommand{\SLF}{b}
\newcommand{\Nbdy}{{N_{\mathrm{bdy}}}}
\newcommand{\npred}{{q}}
\newcommand{\xprime}[2]
{
\ifthenelse{\equal{#1}{\string 1}}
{{\frac{d}{d t}x\left(#2\right)}}
{{\frac{d^{#1}}{dt^{#1}}x\left({#2}\right)}}%
}
\renewcommand{\xprime}[2]
{
\ifthenelse{\equal{#1}{\string 1}}
{{x'\left(#2\right)}}
{{x^{[#1]}\left({#2}\right)}}%
}
\newcommand{\bernhard}[1]{{{\textcolor{magenta}{#1}}}}
\newcommand{\peter}[1]{{\textcolor{olive}{#1}}}
%comment by Peter
\newcommandx{\pcomment}[2][1=]{\todo[linecolor=red,backgroundcolor=red!25,bordercolor=red,#1]{#2}}
%comment by Bernhard
\newcommandx{\bcomment}[2][1=]{\todo[linecolor=olive,backgroundcolor=olive!25,bordercolor=olive,#1]{#2}}
%comment by Alessio
\newcommandx{\acomment}[2][1=]{\todo[linecolor=blue,backgroundcolor=blue!25,bordercolor=blue,#1]{#2}}
\providecommand{\keywords}[1]{\textbf{\textit{Keywords:}} #1}
\providecommand{\msc}[1]{\textbf{\textit{2020 MSC:}} #1}

\newif\ifpreprint
 \preprinttrue         %set \preprinttrue to enable version for arxiv (including example of higher order shape deriv for PDE)

\title{
	Homotopy methods for higher order shape optimization: \\A globalized shape-Newton method and Pareto-front tracing
}

\author[1]{A. Cesarano}
\author[2,3]{B. Endtmayer}
\author[1]{P. Gangl}

\affil[1]{Johann Radon Institute for Computational and Applied Mathematics, Altenbergerstr. 69, A-4040 Linz, Austria}	
\affil[2]{Leibniz Universit\"at Hannover,
	Institut f\"ur Angewandte Mathematik,
	AG Wissenschaftliches Rechnen,
	Welfengarten 1, 30167 Hannover, Germany}
\affil[3]{Cluster of Excellence PhoenixD,
	(Photonics, Optics, and Engineering - Innovation Across Disciplines), 
	Leibniz Universit\"at Hannover, Germany}

\date{}

\maketitle

\begin{abstract}
First order shape optimization methods, in general,
require a large number of iterations until they reach a locally optimal design. While higher order methods can significantly reduce the number of iterations, they exhibit only local convergence properties, necessitating a sufficiently close initial guess. In this work, we present an unregularized shape-Newton method and combine shape optimization with homotopy (or continuation) methods in order to allow for the use of higher order methods even if the initial design is far from a solution. The idea of homotopy methods is to continuously connect the problem of interest with a simpler problem and to follow the corresponding solution path by a predictor-corrector scheme. We use a shape-Newton method as a corrector and arbitrary order shape derivatives for the predictor. Moreover, we apply homotopy methods also to the case of multi-objective shape optimization to efficiently obtain well-distributed points on a Pareto front. Finally, our results are substantiated with a set of numerical experiments.
\end{abstract}

% REQUIRED
\begin{keywords}
continuation,
homotopy,
path following,
nonlinear problem,
shape optimization,
shape Newton,
Pareto tracing,
multi-objective shape optimization
\end{keywords}

% REQUIRED
\begin{msc}
49Q10, %!Optimization of shapes other than minimal surfaces
49M15, %!Newton-type methods
49M41, %!PDE constrained optimization (numerical aspects)
58E17 %! Multiobjective variational problems, Pareto optimality, applications to economics
% 34A12 %~Initial value problems, existence, uniqueness, continuous dependence and continuation of solutions to ordinary differential equations
% 68Q25,%?Analysis of algorithms and problem complexity
% 68R10,%?Graph theory (including graph drawing) in computer science
% 68U05%?Computer graphics; computational geometry (digital and algorithmic aspects)
\end{msc}

%%%%%%%%%%%%%%%%%%%%%%%%%%%%%%%%%%%%%%%%%%%%%%%%%%%%%%%%%%%%
\section{Introduction} \label{sec_intro}
Over the past decades, mathematical and numerical methods for shape optimization have been thoroughly investigated in both the mathematical and the engineering community and have been applied to a wide range of real world problems. We refer to the monographs \cite{b_DEZO_2011a, SZ} for detailed introductions to the main concepts of shape optimization and, exemplarily, mention concrete applications from solid \cite{EtlingHerzog2018} or fluid mechanics \cite{Pinzon2022}, multi-physics problems \cite{Feppon2019Sep, FEPPON2021113638}, electrical engineering \cite{GLLMS2015}, acoustics \cite{SchmidthWadbroBerggren2016} or inverse problems \cite{KovtunenkoOhtsuka2022}.

While most of these works are based on first order shape derivatives, also second order shape derivatives have been investigated from both an analytical and numerical viewpoint \cite{Delfour1991, a_SC_2018a, SchmidtSchulz2023, EtlingHerzogLoayzaWachsmuth2020, sturm2018convergence, GanglEtAlSAMO2020}. A key difficulty when employing Newton methods in shape optimization is the kernel of the shape Hessian operator which includes all interior and tangential shape perturbations. In order to approximately solve the shape Newton system, often regularization techniques are employed which, however, may prevent the desired superlinear or even quadratic convergence unless devised carefully \cite{SchmidtSchulz2023}. As an alternative, the shape Newton equation can be solved for normal perturbations only as proposed in \cite{EtlingHerzogLoayzaWachsmuth2020}. In \cite{sturm2018convergence}, the author proves superlinear convergence of a shape Newton method by employing approximately normal functions.

Another well-known difficulty in the application of Newton's method for (typically non-convex) shape optimization problems is the method's local convergence properties, meaning that Newton's method is only applicable if the initial design is sufficiently close to a solution of the optimization problem. As a simple remedy, one could couple first and second order methods by first running a gradient descent method until a certain residual is reached before switching to a second order method yielding fast convergence. In this approach, however, the threshold that needs to be reached by the first order method may vary from problem to problem, making the method unrobust.

In this paper, we propose to couple shape optimization with homotopy (or continuation) methods \cite{allgower2012numerical, Deuflhard2011} in order to improve convergence properties of shape-Newton methods even if the initial design is far from the sought solution. Homotopy methods have been introduced as a method for solving difficult systems of equations, see e.g. \cite{Malinen2010, Zulehner1988} or \cite{Dunlavy2005, WatsonHaftka1989} for the context of optimization problems.
These methods rely on the idea of smoothly connecting a problem of interest whose solution is sought with a much simpler problem whose solution is known or can easily be computed. 
Given such a connection, called the homotopy map, one can follow the path of solutions to intermediate problems to reach the problem of interest by means of, e.g., a predictor-corrector method.
As shown in \cite{Deuflhard2011, suli2003introduction}, Newton's method just ensures quadratic convergence if the initial guess is in a neighborhood of the solution. The basic idea of the predictor corrector method is to stay in this neighbourhood when following the path.
We will apply this method to the problem of finding stationary shapes with vanishing shape derivative. The employed corrector method will be a shape Newton method and we will use polynomial predictors of first and higher order.

Beside their globalizing effect, there are situations where continuation methods are useful for their own sake as the choice of the homotopy allows to implicitly control the path that is taken during the optimization. Often, optimization problems involve parameters coming from the mathematical model or the algorithm to be used, for which a gradual variation can be beneficial in terms of the obtained solution. As an example we mention the penalization parameter in density-based topology optimization \cite{Houta2023pre}. Moreover, as we will elaborate in Section \ref{sec_pareto}, homotopy methods can be used for finding points on a Pareto front in the case of multi-objective optimization problems \cite{Hillermeier2001JOTA, Hillermeier2001book}, and are also useful for systematically exploring the design space by, e.g., deflation methods \cite{Papadopoulos2021}. Preliminary numerical results related to Section \ref{sec_pareto} have been submitted to a conference proceeding, see \cite{CesaranoGanglSCEE}.

The novelty of this paper is threefold: On the one hand, we present a new way of solving the unregularized Newton system which, to the best of the authors' knowledge, has not been reported in the literature. The main contribution of this work is the coupling of shape optimization with homotopy methods including higher order predictors. Finally, as a by-product, we also obtain an efficient way for tracing Pareto-optimal shapes of multi-objective shape optimization problems.

The remainder of this work is organized as follows: In Section \ref{sec_homotopy}, we present the concept of homotopy methods along with different predictor choices and step size strategies. We briefly discuss concepts of shape optimization in Section \ref{sec_shapeOpti} before presenting the connection of these two concepts in Section \ref{sec_homo_shape}. We comment on its use in the context of multi-objective optimization in \ref{sec_pareto} before presenting a set of numerical experiments in Section \ref{sec_numerics}.

\section{Homotopy methods} \label{sec_homotopy}
The idea of homotopy methods for solving a problem of interest is to establish a smooth connection with a simple problem whose solution is known or can easily be computed, and to follow the path of solutions to the arising intermediate problems until the original problem is reached. For better illustration, consider a function $F : \R^n \rightarrow \R^n$ and let $x^{(0)} \in \R^n$ be given. We look for a root satisfying
\begin{align} \label{eq_Fxzero}
    F(x) = 0.
\end{align}
Often, a solution to problem \eqref{eq_Fxzero} cannot be found from the given initial guess $x^{(0)}$. A homotopy map $H: \R^n \times \R \rightarrow \R^n$ defines a family of equations
\begin{align}
    H(x,t) = 0
\end{align}
parametrized by the homotopy parameter $t \in [0,1]$. The homotopy map should be chosen in such a way that solution to the initial problem $H(x,0) = 0$ can be easily obtained from the initial guess $x^{(0)}$, and that the equation $H(x,1) = 0$ is equivalent to the original equation \eqref{eq_Fxzero}.
We mention one widely used class of homotopy maps which is also considered throughout this article.
\begin{example}[convex homotopy]
    Let $G:\R^n \rightarrow \R^n$ and $x^{(0)}$ given with {$G(x^{(0)}) = 0$}. Then the convex homotopy of $F$ and $G$ can be defined as
    \begin{align} \label{eq_convHomo}
        H(x,t) = tF(x) + (1-t) G(x).
    \end{align}
    Examples for the choice of $G$ include the Keller model $G(x): = F(x)-F(x^{(0)})$ or the choice $G(x) := x-x^{(0)}$.
\end{example}

Given a homotopy map $H$ for the problem of interest $F$ and a suitable initial guess $x^{(0)}$, the solution path can be followed by a predictor-corrector method. The main ingredients are summarized in Algorithm \ref{algo_predCorr}. Here, one iteratively chooses a step size by which the current homotopy value is increased, defines an improved starting guess for the subsequent homotopy value (predictor step) and tries to solve the problem for the updated homotopy value (corrector step). The latter step is often realized by Newton's method, which converges only if the starting guess is sufficiently close to the solution. Thus, the choice of the step size has a large impact on the efficiency of the overall method. If the step $\Delta t^{(k)}$ at iteration $k$ is chosen too large, the corrector step may fail and one may have to repeat the steps for the current homotopy value $t^{(k)}$ with a decreased step size. We comment on the important issue of step size choice in Section \ref{sec_homo_stepsize}. Similarly, a good predictor can yield better starting guesses for the corrector step and thus allow for larger step sizes and shorter total computation times. We comment on this topic in Section \ref{sec_homo_predictor}.

\begin{algorithm}
    \begin{algorithmic}
    \STATE{ Input: $x^{(0)} \in \R^n$, $t^{(0)}=0$, $k=0$.}
    \WHILE{$t^{(k)}<1$}
        \STATE{Choose step size $\Delta t^{(k)}$ and set $t^{(k+1)} = \mathrm{min}(t^{(k)} + \Delta t^{(k)},1)$. }
        \STATE{\textbf{Predictor:} $\tilde x^{(k+1)} = \mathrm{Pred}(x^{(k)},t^{(k)}, t^{(k+1)}-t^{(k)})$}
        \STATE{\textbf{Corrector:} Obtain $x^{(k+1)}$ as solution to $H(x, t^{(k+1)}) = 0$ by iterative method with starting guess $\tilde x^{(k+1)}$.}
        \STATE{$k \leftarrow k+1$}
    \ENDWHILE
    \end{algorithmic}
    \caption{Predictor-corrector algorithm \label{alg: predictor-corrector-algorithm}}
    \label{algo_predCorr}
\end{algorithm}

\begin{remark}[Homotopy methods for optimization] \label{rem_hom_opti}
    The same principles can also be adapted to the solution of optimization problems rather than nonlinear equations. In that setting, $H(\cdot, t)$ would represent a family of optimization problems and the corrector step of Algorithm \ref{algo_predCorr} would be replaced by an iterative optimization algorithm, see also \cite{Dunlavy2005}.
\end{remark}

\begin{remark}
    Throughout this paper, we assume that the solution path can be parametrized by the homotopy parameter $t$ and that no turning points exist. Moreover, we assume the absence of bifurcations. Then, for a given $t \in [0,1]$, there exists a unique solution to the equation $H(\cdot, t)=0$ on the path of interest, which we denote by $x(t)$.
\end{remark}

\subsection{On the choice of predictors} \label{sec_homo_predictor}
\begin{figure}
    \centering
	\includegraphics[width=.5\textwidth, trim=100 375 20 20, clip]{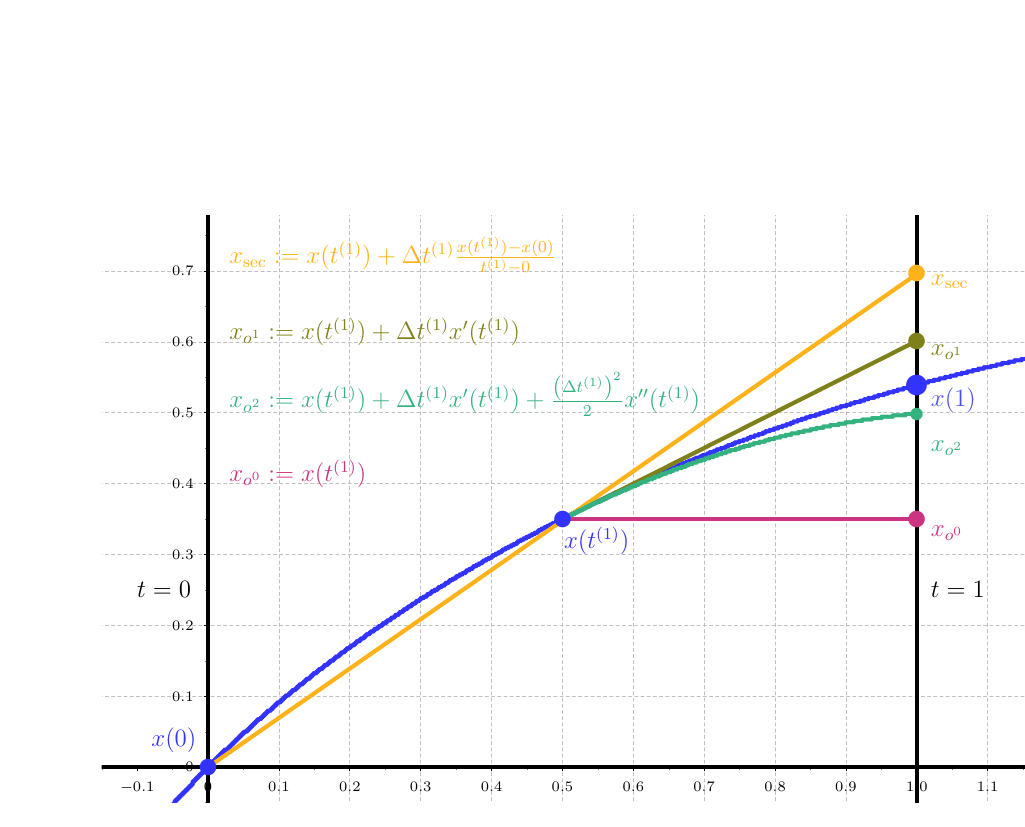}
	\caption{Visualization of a solution curve $(x(t),t)$ (blue) along with different ways of predicting an initial guess starting from $t^{(1)}=0.5$ with $\Delta t^{(1)}=0.5$: Step (magenta), secant (orange), first order (olive), second order (green). Illustration is for problem \eqref{eq_convHomo} with $F(x):=x^5+x-e^{-x}$ and $G(x):=x$.
	}
	\label{fig_predictors}
\end{figure}

Recall that, for $t \in [0,1]$ fixed, $x(t)$ denotes a solution to the homotopy equation for $t$, i.e., it holds
\begin{align}\label{eq_H_xt_t}
    H(x(t), t)=0.
\end{align}
In this section, we assume $x(t)$ to be sufficiently smooth such that all arising derivatives exist. This can be inferred from a sufficient smoothness of $H$ by the implicit function theorem.
For more information on predictor methods for homotopy, we refer to \cite{allgower2012numerical}.
In what follows, we discuss different predictor choices which are also illustrated in Figure \ref{fig_predictors}.

\subsubsection{Zero order predictors}
The most simple predictor is the identity, i.e
\begin{align} \label{eq: pred: do nothing}
\mathrm{Pred}^{(0)}(x^{(k)}, t^{(k)}, \Delta t^{(k)}) = x^{(k)}.
\end{align}
This is an zero order predictor. In some literature, this is also called classical continuation; see \cite{Deuflhard2011}.
Another zero order predictor is the so called secant predictor. In contrast to \eqref{eq: pred: do nothing} this predictor also depends on the solution of the previous step, i.e.
\begin{align} \label{eq: pred: secant}
\mathrm{Pred}^{\textrm{sec}}(x^{(k)}, t^{(k)}, x^{(k-1)}, t^{(k-1)}, \Delta t^{(k)}) = x^{(k)}+\Delta t^{(k)} \frac{x^{(k)}-x^{(k-1)}}{t^{(k)}-t^{(k-1)}}.
\end{align}
For the first step any other predictor can be used. We mention that this predictor does not fit into the class of predictors which are used in Algorithm~\ref{alg: predictor-corrector-algorithm}.
\subsubsection{First order predictors}
Differentiating \eqref{eq_H_xt_t} with respect to $t$ yields
\begin{align} \label{eq_dtH}
    H_x(x(t), t)x'(t) + H_t(x(t),t) = 0,
\end{align}
where $H_x: \R^n\times \R \rightarrow \R^{n\times n}$ denotes the Jacobian of $H(\cdot, t)$ and $H_t: \R^n\times \R \rightarrow \R^n$ the partial derivative of the homotopy map with respect to the homotopy parameter~$t$. Assuming invertibility of $H_x(x(t),t)$, this defines the first order derivative
\begin{align} \label{eq_Davidenko}
    x'(t) =- \left( H_x(x(t), t) \right)^{-1} H_t(x(t),t) =:  g_{\mathrm{Dav}}(x(t), t)
\end{align}
which gives rise to a first order predictor
\begin{align} \label{eq: pred: first order predictor}
    \mathrm{Pred}^{(1)}(x^{(k)}, t^{(k)}, \Delta t^{(k)}) = x^{(k)} + \Delta t^{(k)} g_{\mathrm{Dav}}(x^{(k)}, t^{(k)}).
\end{align}
This is also known as tangent predictor our tangent continuation as explained in \cite{Deuflhard2011}
\begin{remark}
    Instead of the predictor-corrector scheme described above, one could also consider the related initial value problem given by \eqref{eq_Davidenko} together with the initial condition $x(0) = x^{(0)}$ and apply numerical techniques for ordinary differential equations, see \cite{BoltenDoganayGottschalkKlamroth2021}. Equation \eqref{eq_Davidenko} is called the Davidenko differential equation.
\end{remark}

\subsubsection{Higher order predictors}
As it was mentioned above, a better predictor can speed up the overall computation as it may allow for larger step sizes and reduce the number of iterations needed in the corrector step. For smooth enough problems one may thus approximate the unknown solution $x(t+\Delta t)$ by higher order Taylor polynomials involving higher order derivatives of the mapping $t \mapsto x(t)$, see Fig. \ref{fig_predictors} for an illustration. These higher order derivatives $\xprime{n}{t}$, $n\geq 2$ can be obtained by further differentiating \eqref{eq_dtH}.
Note the difference between the employed notation $x^{(k)}=x(t^{(k)})$ and the $n$-th derivative $x^{[n]}(t)$ of the mapping $t \mapsto x(t)$.

For that purpose, let $X^{(k)}(t):= \left\lbrace\xprime{k}{t},\xprime{k-1}{t},\ldots ,\xprime{2}{t}, \xprime{1}{t} \right\rbrace$
denote the set of all derivatives of $t \mapsto x(t)$ up to order $k$. Moreover, let $\alpha_k:=\{(l,m) \in (\mathbb{N}\cup \{0\})^2, l+m\leq k \}$ denote all multi-indices of order $k$ or less and $H_{\alpha_k}:=\left\{\frac{\partial^l}{\partial x^l} \frac{\partial^m}{\partial t^m}H: (l,m) \in \alpha_k \right\}$ the set of all derivatives of $H$ corresponding to $\alpha_k$.
\begin{theorem} \label{thm_recursion_homo}
	Let $H(x(t),t)$ be sufficiently smooth and assume that $H_x(x(t),t)$ is invertible for all $t\in [0,1]$. Then the $n$-th derivative of $x(t)$, $\xprime{n}{t}$, is given as the solution to the system
	\begin{align} \label{eq_def_dtnx}
	H_x(x(t),t) \xprime{n}{t} = \SLF^{(n)}\left(X^{(n-1)}(t),H_{\alpha_n}(x(t),t)\right),
	\end{align}
	where $\SLF^{(n)}(\cdot, \cdot)$ is a semilinear form which satisfies the recurrence relation
	\begin{align*}
		\SLF^{(n+1)}\left(X^{(n)}(t),H_{\alpha_{n+1}}(x(t),t)\right):=&\SLF^{(n)}_X\left(X^{(n-1)}(t),H_{\alpha_n}(x(t),t)\right)X^{(n-1)}_t(t)\\
		+&\SLF^{(n)}\left(X^{(n-1)}(t),H_{\alpha_n x}(x(t),t)\xprime{1}{t}+H_{\alpha_n t}(x(t),t)\right)\\
		-&H_{xx}(x(t),t)\left[\xprime{n}{t},\xprime{1}{t} \right]-H_{xt}(x(t),t)\xprime{n}{t}.
	\end{align*}
	with $b^{(1)}(X^{(0)}(t), H_{\alpha_1}(x(t), t)) = - H_t(x(t), t)$.

	\begin{proof}
		We prove this statement by induction. The base case $n=1$ was already derived in \eqref{eq_Davidenko}.
		Now let us assume the statement is already valid for $n$, i.e., we assume
		\begin{equation*}
		H_x(x(t),t)\xprime{n}{t}=\SLF^{(n)}\left(X^{(n-1)}(t),H_{\alpha_n}(x(t),t)\right).
		\end{equation*}
		Again differentiating with respect to $t$ yields
		\begin{align*}
			H_x(x(t),t)\xprime{n+1}{t}&+H_{xx}(x(t),t) \left[\xprime{n}{t},\xprime{1}{t} \right]+H_{xt}(x(t),t)\xprime{n}{t}\\
			=&\SLF^{(n)}_X\left(X^{(n-1)}(t),H_{\alpha_n}(x(t),t)\right)X^{(n-1)}_t(t) \\
			&+\SLF^{(n)}\left(X^{(n-1)}(t),H_{\alpha_n x}(x(t),t)\xprime{1}{t}+H_{\alpha_n t}(x(t),t)\right),
		\end{align*}
		which, after rearranging, yields
		\begin{equation*}
		H_x(x(t),t)\xprime{n+1}{t}=\SLF^{(n+1)}\left(X^{(n)}(t),H_{\alpha_n}(x(t),t)\right).
		\end{equation*}
	\end{proof}	
\end{theorem}

\begin{remark}
    Of course there is a trade-off between the reduction of computational costs by the use of additional predictor terms and the cost of computing these terms. We note, however, that the system matrix for computing the $n$-th derivative of $t \mapsto x(t)$ in \eqref{eq_def_dtnx} is the same for all $n$. Thus, if the system matrix can be factorized, the factorization has to be done only once and additional derivatives can be computed at the cost of assembling the right hand sides of \eqref{eq_def_dtnx}.
\end{remark}

Given the first $\npred$ derivatives of $t \mapsto x(t)$, an $\npred$-th order predictor at $x^{(k)}=x(t^{(k)})$ can be defined as
\begin{align} \label{eq: pred: higher order}
    \mathrm{Pred}^{(\npred)}(x^{(k)}, t^{(k)}, \Delta t^{(k)}) := x^{(k)} + \sum_{i=1}^\npred \frac{1}{i!}\left( \Delta t^{(k)} \right)^i \xprime{i}{t^{(k)}}.
\end{align}

\begin{corollary} \label{cor_predictors} In addition to the system defining the first order derivative of $t\mapsto x(t)$ stated in \eqref{eq_dtH} which corresponds to \eqref{eq_def_dtnx} for $n=1$, we also give the problems defining the second and third order derivatives:
    \begin{align} \label{eq_dttH}
        H_x[ x''(t)] =& -  H_{xx}[x'(t), x'(t)] - 2 H_{xt}[x'(t)] - H_{tt},\\
    \begin{split}
    \label{eq_dtttH}
        H_x [x'''(t)] =& - \Bigg[H_{xxx}[x'(t), x'(t), x'(t)] +3H_{xxt}[x'(t), x'(t)] \\
        &+ 3 H_{xx}[x''(t), x'(t)] +3  H_{xt}[x''(t)] + 3 H_{xtt}[x'(t)]+H_{ttt}\Bigg].
    \end{split}
    \end{align}
    Note that the arguments $(x(t), t)$ were omitted here for brevity.
\end{corollary}

\subsection{Step size choice} \label{sec_homo_stepsize}
The choice of the step sizes $\Delta t^{(k)}$ in Algorithm \ref{algo_predCorr} has a big impact on the overall efficiency of the method. While larger step sizes reduce the number of necessary homotopy steps to reach the solution of the original problem at $t=1$, the quality of the starting guesses for the corrector phase is decreased and the corrector step may be unsuccessful.

A straightforward way to define step sizes would be to equidistant steps, i.e., to choose $\Delta t^{(k)} = 1/m$ for a fixed number $m$ of steps. This choice, however, can be inefficient when small steps are necessary for a certain part of the homotopy path, but larger steps could be made in other parts. Thus, we propose three different ways of adapting the step sizes along the way:
\subsubsection{Fixed adaptation of step sizes} \label{sec_fixed_adap}
We start with $t^{(0)}=0$ and an initial step size $\Delta t^{(0)}>0$. If, for $k\geq 1$, the combination of predictor and corrector step to get from homotopy value $t^{(k-1)}$ to $t^{(k)}$ by the increment $\Delta t^{(k-1)}$ was successful, we increase the step size for the next step by a constant factor $\overline \gamma>1$, i.e., $\Delta t^{(k)} = \overline \gamma \Delta t^{(k-1)}$. If, on the other hand, the corrector step is not successful, we reduce the step size by a constant factor $\underline \gamma<1$ and repeat the corrector step starting from the adapted predicted value.

\subsubsection{Agile adaptation of step sizes based on predictors} \label{sec_step_agile}
While the strategy of Section \ref{sec_fixed_adap} is capable of adapting the step size to the problem, this adaptation happens only gradually. In particular, sudden changes in the desirable step size or the detection of a suitable initial step size $\Delta t^{(0)}$ are not possible in this way. For that reason, we propose to choose the step size depending on the estimated distance between the solution at $t = t^{(k+1)}$ and its prediction of a given order. More precisely, given a predictor of order $\npred$ via \eqref{eq: pred: higher order}, by Taylor expansion it holds
\begin{align*}
    \left \lvert x \left(t^{(k+1)} \right) - \mathrm{Pred}^{(\npred)} \left(x^{(k)}, t^{(k)}, \Delta t^{(k)} \right) \right \rvert = \frac{(\Delta t^{(k)})^{\npred+1}}{(\npred+1)!}  \left \lvert \xprime{\npred+1}{t^{(k)}} \right \rvert + \mathcal O((\Delta t^{(k)})^{\npred+2})
\end{align*}
as $\Delta t^{(k)} \rightarrow 0$. The idea now is, at the cost of computing one additional derivative, to choose $\Delta t^{(k)}$ such that the leading term on the right hand side becomes a constant $\alpha>0$ (which is in particular independent of the curve $t \mapsto x(t)$), i.e., we choose
\begin{align} \label{eq_agile_step_size}
    \Delta t^{(k)} = \left( (\npred+1)! \, \alpha \right)^{1/(\npred+1)} \left \lvert \xprime{\npred+1}{t^{(k)}} \right \rvert^{-1/(\npred+1)}.
\end{align}
By the constant $\alpha$, the user can determine the level of conservativeness. Smaller values of $\alpha$ will increase the success probability of the corrector step, but will increase the total number of homotopy steps.

\subsubsection{Agile adaptation with adaptive parameter choice} \label{sec_step_agileadap}
While the strategy of Section \ref{sec_step_agile} allows for an automatic detection of the order of magnitude of a step size, the best possible choice of the parameter $\alpha$ is again problem-dependent. For that reason, we suggest to combine the strategy of Section \ref{sec_step_agile} with the parameter adaptation described in Section \ref{sec_fixed_adap}, i.e. to multiply $\alpha$ by constants $\underline \alpha<1$ or $\overline \alpha>1$ in case of failure or success of the corrector step, respectively.

\section{Shape optimization} \label{sec_shapeOpti}
We briefly review the main concepts needed for first and second order shape optimization algorithms. For a more thorough introduction, we refer the reader to the monographs \cite{b_DEZO_2011a, SZ} or the review papers \cite{AllaireDapogny2021, GanglEtAlSAMO2020}.

Numerical shape optimization algorithms are based on the concept of shape derivatives, i.e., sensitivities of a shape-dependent cost function with respect to a variation of the domain represented by the action of deformation vector fields. To be more precise, let $\Dsf \subset \R^d$ denote a hold-all domain, let $\mathcal A \subset \mathcal P(\Dsf)$ be a set of admissible shapes and consider a shape function $\mathcal J: \mathcal A \rightarrow \R$. For a shape $\Omega \in \mathcal A$, a smooth vector field $\V: \R^d \rightarrow \R^d$ and a parameter $s>0$, let $\Omega_s := (\mathrm{id} +s \V)(\Omega)$ denote the corresponding perturbed domain. The first order shape derivative of the shape function $\mathcal J$ is defined as
\begin{align} \label{eq_def_dJ}
    d \mathcal J(\Omega)(\V) := \left.\left(\frac{d}{ds} \mathcal J(\Omega_s) \right) \right\rvert_{s=0} = \underset{s \rightarrow 0}{\mbox{lim }} \frac{\mathcal J(\Omega_s) - \mathcal J(\Omega)}{s}
\end{align}
if this limit exists and the mapping $\V \mapsto d\mathcal J(\Omega)(\V)$ is linear.
Similarly, we define the symmetric second order shape derivative by considering smooth vector fiels $\V, \W : \R^d \rightarrow \R^d$, two parameters $s, r >0$ and the corresponding perturbed domain $\Omega_{s,r}:= (\mathrm{id}+s \V + r \W)(\Omega)$. The symmetric second order shape derivative is defined as
\begin{align} \label{eq_def_d2J}
    d^2 \mathcal J(\Omega)[\V, \W] := \left.\left(\frac{d^2}{dsdr} \mathcal J(\Omega_{s, r}) \right) \right\rvert_{s=0, r=0}.
\end{align}
if it exists and the mapping $(\V,\W)\mapsto d^2 \mathcal J(\Omega)[\V, \W]$ is bilinear.
We mention that there also exists a non-symmetric second order shape derivative which is obtained by repeated shape differentiation, i.e.,
\begin{align} \label{eq_shapeHess_nonsymm}
    D^2 \mathcal J(\Omega)[\V][ \W] := \underset{r \rightarrow 0}{\mbox{lim }} \frac{1}{r} \bigg( d \mathcal J( (\mathrm{id}+r \W)(\Omega))(\V) - d \mathcal J(\Omega)(\V) \bigg)
\end{align}
with the relation \cite{b_DEZO_2011a, sturm2018convergence, SchmidtSchulz2023, Simon1989}
\begin{align*}
    D^2 \mathcal J(\Omega)[\V][\W] = d^2 \mathcal J(\Omega)[\V, \W] + d \mathcal J(\Omega)(\partial \V \W).
\end{align*}

\subsection{First order methods}
The main idea of most first order shape optimization algorithms consists in repeatedly finding a descent direction $\V$ such that $d\mathcal J(\Omega)(\V)<0$ and updating the current domain by the action of this vector field to get $\tilde{\Omega} = (\mathrm{id} + s \V)(\Omega)$ as the new iterate. Here, the parameter $s>0$ is chosen sufficiently small in order to guarantee a descent, $\mathcal J(\tilde{\Omega})< \mathcal J(\Omega)$. A descent direction can be obtained by solving an auxiliary boundary value problem to find $\V \in H$ satisfying
\begin{align} \label{eq_desc_direction}
    b(\V, \W) = - d \mathcal J(\Omega)(\W)
\end{align}
for all $\W \in H$ for a given vector-valued Hilbert space $H$ and a positive definite bilinear form $b(\cdot, \cdot) : H \times H \rightarrow \R$. Possible choices for $b(\cdot, \cdot)$ include a vector Laplacian or elasticity problem on (a subdomain of) $H^1(\Dsf)^d$, the latter being advantageous in terms of mesh quality \cite{SchulzSiebenbornWelker2016}. Other means to ensure mesh quality of deformed finite element meshes include the use of nearly conformal mappings \cite{a_IGSTWE_2018a} or solving nonlinear versions of \eqref{eq_desc_direction} \cite{Onyshkevych2021, Mueller2021} or even employing $W^{1, \infty}$ descent directions \cite{DeckelnickHerbertHinze2022}.

\subsection{Second order methods}
Replacing the bilinear form $b(\cdot, \cdot)$ in \eqref{eq_desc_direction} by the (symmetric) second order shape derivative defined in \eqref{eq_def_d2J} for the problem at hand and the current domain $\Omega$ yields the shape-Newton system to find $\V \in H$ satisfying
\begin{align}\label{eq_shapeNewton}
    d^2\mathcal J(\Omega)(\V, \W) = - d \mathcal J(\Omega)(\W)
\end{align}
for all $\W \in H$ for a suitable Hilbert space $H$. The operator on the left hand side of \eqref{eq_shapeNewton}, however, is not invertible. In fact, its kernel contains all vector fields whose normal component on the boundary $\partial \Omega$ vanishes \cite{SchmidtSchulz2023, NovruziPierre2002, sturm2018convergence}. This includes vector fields corresponding to deformations which only move interior points as well as tangential movements. In order to obtain an approximate shape-Newton descent direction, problem \eqref{eq_shapeNewton} is often modified either by regularizing the operator on the full domain \cite{a_SC_2018a} or only on the boundary \cite{GanglEtAlSAMO2020}.

\subsection{Numerical shape-Newton method} \label{sec_shapeNewton}
In this work, we propose to solve the finite-element discretized version of \eqref{eq_shapeNewton} without regularization, by filtering out the full kernel of the shape Hessian. To that end, consider a given initial domain $\Omega \subset \R^d$ which is discretized by a conforming simplicial finite element mesh $\mathcal T_h$. We introduce the space of piecewise linear and globally continuous functions on $\mathcal T_h$,
\begin{align} \label{eq_def_Vh}
    V_h := \{ v \in H^1(\Omega)\cap C(\bar{\Omega}): v|_T \in P^1(T) \text{{ for all }} T \in \mathcal T_h\}
\end{align}
where $P^1(T)$ denotes the space of linear functions on a simplex $T$ (i.e. a triangle if $d=2$ or a tetrahedron if $d=3$). We denote the number of vertices of $\mathcal T_h$ by $N$ and the number of vertices that lie on the boundary of $\Omega$ by $\Nbdy$. The linear hat basis function of $V_h$ corresponding to vertex $x_k$ is denoted by $\varphi_k$, $k=1, \dots, N$. Without loss of generality, we assume the basis functions to be numbered in such a way that the first $\Nbdy$ basis functions correspond to the vertices on the boundary.

Our goal is now to find a solution $\V_h \in V_h^d$ to the discretized version of \eqref{eq_shapeNewton}. We deal with the kernel of this problem as follows:
\begin{itemize}
 \item Interior deformations are removed by choosing the Hilbert space $H$ in \eqref{eq_shapeNewton} to be defined solely on the surface, $H := (V_h)^d \cap H^1(\partial \Omega)^d$ with $\mathrm{dim}(H) = d \Nbdy$ and basis functions (for $k=1, \dots, \Nbdy$): $\Phi_{2k-1} = (\tilde \varphi_k,0)^\top$, $\Phi_{2k} = (0, \tilde \varphi_k)^\top$ in the case $d=2$ and $\Phi_{3k-2} = (\tilde \varphi_k,0,0)^\top$, $\Phi_{3k-1} = (0, \tilde \varphi_k,0)^\top$ and $\Phi_{3k} = (0, 0, \tilde \varphi_k)^\top$ in the case $d=3$ where $\tilde \varphi_k := \varphi_k|_{\partial \Omega}$.
 \item In order to filter out tangential movements on the boundary, we impose the condition
 $\int_{\partial \Omega} (\V_h \cdot \tau)(x) \varphi_k(x) \; \mbox{d}S = 0$, $k=1, \dots, \Nbdy$, which implies the pointwise condition $(\V_h \cdot \tau)(x_k) = 0$ for all mesh vertices $x_k$ on the boundary of $\Omega$. Here, the unit tangential vector $\tau$ on $\partial \Omega$ at a boundary vertex $x_k$ is obtained by averaging of the adjacent edges/faces.
\end{itemize}
Summarizing, we solve the problem to find $\V_h = \sum_{j=1}^{d\Nbdy} \V_j \Phi_j \in H$ such that $\underline V := (\V_j)_{j=1}^{d\Nbdy}$ together with $\underline \xi \in \R^\Nbdy$ satisfies
\begin{align} \label{eq_shapeNewton_Vn}
    \left( \begin{array}{cc}
        \left[d^2\mathcal J(\Omega)(\Phi_j, \Phi_i) \right]_{i,j=1}^{d \Nbdy} & B_{\tau, \Omega} \\
        B_{\tau, \Omega}^\top & 0
    \end{array} \right) \left( \begin{array}{c}
                                 \underline V \\
                                 \underline \xi
                               \end{array} \right) =
    \left( \begin{array}{cc}
            - \left[ d\mathcal J(\Omega)(\Phi_i) \right]_{i=1}^{d\Nbdy} \\ 0
           \end{array} \right).
\end{align}
with the matrix $B_{\tau, \Omega}$ defined by
\begin{align} \label{eq_def_BtauOmega}
    \big[ B_{\tau, \Omega} \big]_{i,k} = \int_{\partial \Omega}(\Phi_i \cdot \tau)\varphi_k \mbox{d}S, \quad i=1, \dots, d\Nbdy, \; k = 1, \dots \Nbdy.
\end{align}
The solution of \eqref{eq_shapeNewton_Vn} yields a vector field $\V_h \in H$ that is normal to $\partial \Omega$. In a next step, we compute a smooth extension $\hat \V_h$ of $\V_h$ by solving the Dirichlet problem of linear elasticity type to find $\hat \V_h \in \{\mathbf v \in V_h^d: \mathbf v|_{\partial \Omega} = \V_h\}$ such that
\begin{align}\label{eq_extension}
   \int_\Omega 2 \mu \varepsilon(\hat \V_h) : \nabla \mathbf v + \lambda \mbox{div}(\hat \V_h) \mbox{div}(\mathbf v)  \; \mbox{d}x = 0 \quad \mbox{ for all } \mathbf v \in V_h^d
\end{align}
with the linearized strain tensor $\varepsilon(\mathbf u) = \frac{1}{2}(\nabla \mathbf u + \nabla \mathbf u^\top)$ and Lam\'e parameters $\mu, \lambda$.

The proposed procedure is summarized in Algorithm \ref{algo_shapeNewton}. As a stopping criterion, we set a tolerance to the $L^2(\partial \Omega)$ norm of the Newton update $\V_h$. We mention that our algorithm is related to the method proposed in \cite{EtlingHerzogLoayzaWachsmuth2020}.
\begin{algorithm}
    \begin{algorithmic}
        \STATE{Input: $\Omega^{(0)} \in \mathcal A$, $\mathrm{tol}$, $iter_{\mathrm{max}}$, $j=0$.}
        \WHILE{$j<iter_{\mathrm{max}}$}
            \STATE{Assemble system matrix and right hand side $d\mathcal J(\Omega^{(j)})$ of \eqref{eq_shapeNewton_Vn}}
            \STATE{Solve system \eqref{eq_shapeNewton_Vn} to obtain $\V_h \in H$}
            \STATE{Extend $\V_h$ to $\hat \V_h \in V_h$ by \eqref{eq_extension}}
            \STATE{Update shape: $\Omega^{(j+1)} \leftarrow (\text{id}+ \hat \V_h)(\Omega^{(j)})$, $j\leftarrow j+1$}
            \IF{$\|\V_h\|_{L^2(\partial \Omega)} < \mathrm{tol}$}
                \STATE \textbf{break}
            \ENDIF
        \ENDWHILE
    \end{algorithmic}
    \caption{Shape-Newton algorithm without regularization}
    \label{algo_shapeNewton}
\end{algorithm}

\subsection{Arbitrary order shape differentiation} \label{sec_shapeDiff_higher}
In order to treat higher order predictors, we will also deal with higher order shape derivatives.
\begin{definition}
Let a domain $\Omega$, $k$ smooth vector fields $\V^{(1)}, \dots, \V^{(k)}$, and $s \in \R^k$ with $s_i \geq 0$, $i=1,\dots, k$ be given. Let the corresponding perturbed domain
$\Omega_s := \left(\mathrm{id}+ \sum_{i=1}^k s_i \V^{(i)}\right)(\Omega).$
The symmetric $k$-th order shape derivative of a shape function $\mathcal J(\Omega)$ at $\Omega$ in the directions of $\V^{(1)}, \dots, \V^{(k)}$ is defined as
\begin{align}
    d^k \mathcal J(\Omega)[\V^{(1)}, \dots, \V^{(k)}] := \left. \frac{d^k}{ds_1 \cdots ds_k} \mathcal J(\Omega_s) \right\rvert_{s_1=\dots=s_k=0}.
\end{align}
\end{definition}
We illustrate how a higher order shape derivative can be computed in practice for a simple shape optimization problem.

\begin{example} \label{ex_intf_highershapeder}
    Let $f: \R^d \rightarrow \R$ smooth and consider the cost function
    \begin{align}
        \mathcal J(\Omega) = \int_\Omega f(x) \; \mbox{d}x
    \end{align}
    for smooth domains $\Omega$. Given $s \in \R^k$, $s \geq 0$ componentwise and the smooth vector fields $\V^{(1)}, \dots, \V^{(k)} : \R^d \rightarrow \R^d$, we define the transformation $T_s(x) = (\mathrm{id}+ \sum_{i=1}^k s_i \V^{(i)})(x)$ with Jacobian $F_s(x) := \partial T_s(x) = I + \sum_{i=1}^k s_i \partial \V^{(i)}(x)$. The $k$-th order shape derivative
    can now be obtained by transforming the perturbed domain $\Omega_s$ back to the original domain $\Omega$ and repeatedly differentiating the arising integrand with respect to $T_s$ in the direction of $\frac{d T_s}{ds_i} = \V^{(i)}(x)$ and with respect to $F_s$ in the direction of $\frac{d F_s}{ds_i} =\partial \V^{(i)}(x)$ for $i=k,k-1, \dots, 1$ at $s=0$, i.e.,
    \begin{align*}
        d^k \mathcal J(\Omega)[\V^{(1)}, \dots, \V^{(k)}]
        =& \left. \frac{d^k}{ds_1 \cdots ds_k} \int_{\Omega} (f\circ T_s)(x) \mbox{det}(F_s) \; \mbox{d}x \right\rvert_{s=0}\\
        \begin{split}
        =&  \frac{d^k}{ds_1 \cdots ds_{k-1}} \int_{\Omega} \left( (\nabla f)\circ T_s(x) \cdot \V^{(k)}(x) \mbox{det}(F_s)  \right. \\
        &\left. \left. +  (f\circ T_s)(x)  \frac{d }{d F_s} \left(\mbox{det}(F_s)\right)\partial \V^{(k)}(x)  \right)\; \mbox{d}x \right\rvert_{s=0} = \dots
        \end{split}
    \end{align*}
    Here it is important to note that the introduced vector fields $\V^{(i)}$ are functions of the unperturbed variable $x$ and not of the perturbed variable $T_s(x)$, thus they do not yield extra terms when being differentiated with respect to $T_s$. This distinguishes the symmetric higher order shape derivative computed here from the quantity obtained by repeated shape differentiation as in \eqref{eq_shapeHess_nonsymm}. We refer the reader to \cite[Sec. 3.2]{GanglEtAlSAMO2020} for a more detailed description and an implementation in the finite element software package NGSolve \cite{Schoeberl2014}. In that framework, using symbolic differentiation of the integrands, arbitrary order symmetric shape derivatives can be computed by repeating the procedure described above. This framework will also be exploited in the numerical experiments of Section \ref{sec_numerics}.
\end{example}

\ifpreprint
\begin{example}
    As we will make use of higher order shape derivatives also in the case of PDE-constrained shape optimization, we also illustrate their computation here. We consider a problem of the form
    \begin{align} \label{eq_exampePDE}
        \underset{\Omega \in \mathcal A}{\mbox{min }} J(\Omega, u) \quad \mbox{ s.t. } \quad u \in V: e(\Omega, u) = 0 \mbox{ in } V^*
    \end{align}
    and introduce the corresponding Lagrangian
    \begin{align}
        L(\Omega, u, p) :=  J(\Omega, u) + e(\Omega, u)(p).
    \end{align}
    Assuming unique solvability of the PDE constraint for each admissible $\Omega$ and denoting the corresponding solution with $u(\Omega)$, we can also introduce the reduced cost function
    \begin{align}
        \mathcal J(\Omega) := J(\Omega, u(\Omega)).
    \end{align}
    It is well-known \cite{a_ST_2015a} that the first order shape derivative of \eqref{eq_exampePDE} is given by the partial derivative with respect to the shape of the Lagrangian when the corresponding state and adjoint are inserted, i.e.,
    \begin{align}
        d \mathcal J(\Omega)[\V] = \partial_\Omega L(\Omega, u(\Omega), p(\Omega))[\V]
    \end{align}
    where $u(\Omega)$ and $p(\Omega)$, respectively, solve
    \begin{align} \label{eq_dpLag}
        \mbox{Find }u \in V: \partial_p L(\Omega, u, 0)[\hat p] =& 0 \quad \mbox{ for all } \hat p \mbox{ in } V, \\
        \mbox{Find }p \in V: \partial_u L(\Omega, u, p)[\hat u] =& 0 \quad \mbox{ for all } \hat u \mbox{ in }V. \label{eq_duLag}
    \end{align}
    Here, the partial shape derivative $\partial_\Omega$ is to be understood as
    \begin{align*}
        \partial_\Omega L(\Omega, w, q)[\V] := \left. \frac{d}{ds} L(T_s(\Omega), \Psi_s(w), \Psi_s(q))\right\rvert_{s=0}
    \end{align*}
    with $T_s(x) = (\mathrm{id}+s\V)(x)$ where the transformation $\Psi_s$ depends on the functional setting (e.g., $\Psi_s(w) = w \circ T_s^{-1}$ if $w \in H^1(\Omega)$), see also \cite[Sec. 4]{GanglEtAlSAMO2020}. The second order shape derivative is formally obtained by an application of the chain rule,
    \begin{align} \label{eq_d2J_examplePDE}
        d^2 \mathcal J(\Omega)[\V, \W] = \partial_\Omega \partial_\Omega L [\V, \W] + \partial_u \partial_\Omega L [\V, u'[\W]] + \partial_p \partial_\Omega L [\V, p'[\W]]
    \end{align}
    where the material derivatives $u'[\W] \in V$ and $p'[\W] \in V$ can be obtained by differentiation of \eqref{eq_dpLag}--\eqref{eq_duLag} with respect to the shape $\Omega$ in the direction $\W$,
    \begin{align} \label{eq_uprime_Lag}
        \partial_u \partial_p L(\Omega, u, 0)[\hat p, u'[\W]] =& -\partial_\Omega \partial_p L(\Omega, u, 0)[\hat p, \W]  \quad \mbox{ for all } \hat p \mbox{ in } V, \\
        \begin{split} \label{eq_pprime_Lag}
        \partial_p \partial_u L(\Omega, u, p)[\hat u, p'[\W]] =& -\partial_\Omega \partial_u L(\Omega, u, p)[\hat u, \W] \\
        &- \partial_u \partial_u L(\Omega, u, p)[\hat u, u'[\W]] \quad \mbox{ for all } \hat u \mbox{ in }V.
        \end{split}
    \end{align}
    Here, we omitted the arguments $(\Omega, u(\Omega),p(\Omega))$ in \eqref{eq_d2J_examplePDE} and \eqref{eq_pprime_Lag} and $(\Omega, u(\Omega),0)$ in \eqref{eq_uprime_Lag} and wrote $u'[\W]$, $p'[\W]$ instead of $u'(\Omega)[\W]$, $p'(\Omega)[\W]$ for brevity. In the same manner, the third order shape derivative is obtained as
    \begin{align*}
         \begin{split}
        d^3 \mathcal J&(\Omega)[\V, \W] = \\
        &\partial_\Omega \partial_\Omega \partial_\Omega L [\V, \W, \Z]  + \partial_u \partial_\Omega \partial_\Omega L [\V, \W, u'[\Z]]  + \partial_p \partial_\Omega \partial_\Omega L [\V, \W, p'[\Z]] \\
        &+ \partial_\Omega \partial_u \partial_\Omega L [\V, u'[\W], \Z] + \partial_u \partial_u \partial_\Omega L [\V, u'[\W], u'[\Z]] + \partial_p \partial_u \partial_\Omega L [\V, u'[\W], p'[\Z]] \\
        &+ \partial_\Omega \partial_p \partial_\Omega L [\V, p'[\W], \Z] + \partial_u \partial_p \partial_\Omega L [\V, p'[\W], u'[\Z]] \\
        &+  \partial_u \partial_\Omega L [\V, u''[\W, \Z]] + \partial_p \partial_\Omega L [\V, p''[\W, \Z]]
         \end{split}
    \end{align*}
    where the second order material derivatives $u''[\W, \Z], p''[\W, \Z] \in V$ can be obtained by further differentiating \eqref{eq_uprime_Lag}, \eqref{eq_pprime_Lag}, i.e., by solving
    \begin{align*}
        \begin{split}
        \partial_u \partial_p &L[\hat p, u''[\W, \Z]]
        = -\partial_\Omega \partial_\Omega \partial_p L[\hat p, \W, \Z] -\partial_u \partial_\Omega \partial_p L[\hat p, \W, u'[\Z]] \\
        & \qquad \qquad\qquad \quad \; - \partial_\Omega \partial_u \partial_p L[\hat p, u'[\W], \Z]
        -\partial_u \partial_u \partial_p L[\hat p, u'[\W], u'[\Z]],
        \end{split} \\
        \begin{split}
         \partial_p \partial_u &L[\hat u, p''[\W, \Z]] = - \partial_u \partial_u L[\hat u, u''[\W, \Z]]\\
         &-\partial_\Omega \partial_\Omega \partial_u L[\hat u, \W, \Z]
         -\partial_u \partial_\Omega \partial_u L[\hat u, \W, u'[\Z]]
         -\partial_p \partial_\Omega \partial_u L[\hat u, \W, p'[\Z]] \\
         &-\partial_\Omega \partial_u \partial_u L[\hat u, u'[\W], \Z]
         -\partial_u \partial_u \partial_u L[\hat u, u'[\W], u'[\Z]]
         -\partial_p \partial_u \partial_u L[\hat u, u'[\W], p'[\Z]] \\
         &-\partial_\Omega \partial_p \partial_u L[\hat u, p'[\W], \Z]
         -\partial_u \partial_p \partial_u L[\hat u, p'[\W], u'[\Z]],
%          -\partial_p \partial_p \partial_u L[\hat u, p'[\W], p'[\Z]],
        \end{split}
    \end{align*}
    for all $\hat p$ and $\hat u$ in $V$, respectively. Note that we used the fact that second and higher order derivatives with respect to the third argument of the Lagrangian vanish by definition.
\end{example}
\else
The case of PDE constraints can be treated similarly, see \cite[Sec. 4]{GanglEtAlSAMO2020} for the case $k=2$ which can be extended to higher order.
\fi

% #######################################################
% ################# SEC 4: Homotopy shape ####################
% #######################################################
\section{Homotopy methods for second order shape optimization} \label{sec_homo_shapeOpti}
In this section, we apply the concepts of homotopy methods introduced in Section \ref{sec_homotopy} to the setting of solving shape optimization problems rather than solving systems of nonlinear equations. We will introduce two model problems in Section \ref{sec_homo_shape} before discussing a predictor-corrector scheme in Section \ref{sec_predCorr_shape}.

\subsection{Homotopies for shape optimization} \label{sec_homo_shape}
Here, we discuss possible ways of defining homotopies for shape optimization problems. We begin with an academic unconstrained shape optimization problem in Section \ref{sec_homo_shape_aca} and discuss the practically more relevant case of PDE-constrained shape optimization problems in Section \ref{sec_homo_shape_pde}.

\subsubsection{An academic shape optimization problem} \label{sec_homo_shape_aca}
We begin by considering a shape optimization problem of the form
\begin{align} \label{eq_shapeOpti_aca}
    \underset{\Omega \in \mathcal A}{\mbox{min}}\, \mathcal J_F(\Omega) = \int_\Omega f(x) \; \mbox{d}x
\end{align}
for a given smooth function $f:\R^d \rightarrow \R$. Here, $\mathcal A$ denotes a set of admissible subsets of $\R^d$. It is straightforward to see that the exact solution of \eqref{eq_shapeOpti_aca} is given by $\Omega^* = \{x \in \R^d: f(x)<0\}$. Given sufficiently smooth vector fields $\V, \W$, first and second order shape derivatives are readily obtained as \cite{GanglEtAlSAMO2020}
\begin{align}
    d \mathcal J_F(\Omega)(\V) =& \int_\Omega \nabla f \cdot \V + f \, \mbox{div} \V \; \mbox{d}x, \\
    \nonumber
    d^2 \mathcal J_F(\Omega)(\V, \W) =& \int_\Omega \nabla^2 f \, \V \cdot \W + \nabla f \cdot \W \,  \mathrm{div}\V + \nabla f \cdot \V \, \mathrm{div}\W \\ & \quad \; + f \, \mathrm{div}\V \, \mathrm{div}\W - f \partial \V^\top : \partial \W \; \mbox{d}x.
\end{align}

As noted in Remark \ref{rem_hom_opti}, the concepts of Section \ref{sec_homotopy} can be transfered to the context of optimzation with minor modifications. Thus, we define the convex combination of shape functionals
\begin{align}\label{eq_defH_aca}
    \mathcal H(\Omega, t) :=  t \mathcal J_F(\Omega) + (1-t) \mathcal J_G(\Omega)
\end{align}
and the corresponding family of shape optimization problems
\begin{align} \label{eq_defHopt_aca}
    \underset{\Omega \in \mathcal A}{\mbox{min}}\, \mathcal H(\Omega, t),
\end{align}
where $\mathcal J_G(\Omega)$ is a simple shape functional with known minimizer. In particular, when starting with the initial design $\Omega^{(0)}$, one may choose
\begin{align} \label{eq_defG_aca}
    \mathcal J_G(\Omega) = \int_\Omega \psi[\Omega^{(0)}] \; \mbox{d}x
\end{align}
where $\psi[\Omega^{(0)}]$ is a level set function representing $\Omega^{(0)}$, i.e., $\Omega^{(0)} = \{x \in \R^d: \psi[\Omega^{(0)}](x) < 0\}$. Obviously, $\Omega^{(0)}$ minimizes the functional $\mathcal J_G$ by construction.

\begin{remark}
    We remark that, of course, also other choices of homotopies are possible. In particular one could define the homotopy directly on the optimality condition of problem \eqref{eq_shapeOpti_aca} and apply classical homotopies such as the Keller model to finding roots of $H(\Omega, t):= d \mathcal J_F(\Omega)(\cdot) - (1-t) d\mathcal J_F(\Omega^{(0)})(\cdot)$. While this approach is independent of an artificially chosen auxiliary problem, the resulting homotopy equation is posed in the dual of the space of vector fields. In order to filter out interior or tangential deformations, the equation would have to be restricted to normal vector fields on the boundary. For simplicity, we chose to restrict ourselves to the class of homotopies defined above, but believe that the Keller model deserves further investigation also in the context of shape optimization.
\end{remark}

\subsubsection{A homotopy for PDE-constrained shape optimization} \label{sec_homo_shape_pde}
Similarly to the case of academic shape optimization problems of the type \eqref{eq_shapeOpti_aca}, we can also deal with PDE-constrained shape optimization problems of the form
\begin{align}
    \label{eq_shapeOpti_pde}
    \underset{\Omega \in \mathcal A}{\mbox{min}}\; J_F(\Omega, u) \quad \mbox{s.t. }u \in V: e_F(\Omega, u) = 0 \; \mbox{ in }V^*
\end{align}
where $V$ denotes some Hilbert space, $J_F: \mathcal A \times V \rightarrow \R$ is the cost function and $e_F(\Omega, \cdot): V \rightarrow V^*$ represents the PDE constraint which we assume to have a unique solution for any given $\Omega \in \mathcal A$. In order to define a homotopy, one can introduce an auxiliary PDE-constrained shape optimization problem whose solution can be easily obtained as
\begin{align}
    \label{eq_defG_pde}
    \underset{\Omega \in \mathcal A}{\mbox{min}}\; J_G(\Omega, u) \quad \mbox{s.t. }u \in V: e_G(\Omega, u) = 0 \; \mbox{ in }V^*
\end{align}
Here, $J_G: \mathcal A \times V \rightarrow \R$ is a simple cost function and $e_G(\Omega, \cdot): V \rightarrow V^*$ represents a simple PDE constraint which admits a unique solution for a given $\Omega \in \mathcal A$.

Defining $J_H(\Omega, u, t):= t J_F(\Omega, u) + (1-t) J_G(\Omega,u)$ and $e_H(\Omega, u, t) := t e_F(\Omega, u) + (1-t)e_G(\Omega, u)$ for $t \in [0,1]$, we can now define a family of PDE-constrained shape optimization problems
\begin{align} \label{eq_defH_pde}
    \underset{\Omega \in \mathcal A}{\mbox{min }} J_H(\Omega, u, t) \quad
    \mbox{s.t. } u \in V: e_H(\Omega, u, t) = 0 \mbox{ in } V^*.
\end{align}
Assuming unique solvability of the PDE constraint $e_H(\Omega, u,t)=0$ for given $\Omega \in \mathcal A$ and $t \in [0,1]$ and denoting the unique solution by $u(\Omega, t)$, we can define the reduced cost function $\mathcal H(\Omega, t):= J_H(\Omega, u(\Omega, t), t)$ and rewrite problem \eqref{eq_defH_pde} as the family of shape optimization problems
\begin{align} \label{eq_defH_pde_reduced}
    \underset{\Omega \in \mathcal A}{\mbox{min }} \mathcal H(\Omega, t).
\end{align}

\begin{remark}
    A simple concrete choice for \eqref{eq_defG_pde}, which we also used in the numerical results of Section \ref{sec_numerics}, is given by
   \begin{align} \label{eq_defG_pde_concrete}
        J_G(\Omega, u) =& \int_\Omega u \; \mbox{d}x, \qquad
        e_G(\Omega, u)  = \int_\Omega (u - \psi[\Omega^{(0)}]) (\cdot) \; \mbox{d}x.
   \end{align}
   Noting that the solution $u$ to the PDE constraint is just the $L^2$-projection of $\psi[\Omega^{(0)}]$, we see that $\Omega^{(0)}$ is very close to the solution of \eqref{eq_defG_pde_concrete} in analogy to \eqref{eq_defG_aca}.
\end{remark}

\begin{remark}
    An alternative and possibly more straightforward approach would be to define the reduced cost functions $\mathcal J_F(\Omega)$ and $\mathcal J_G(\Omega)$ of \eqref{eq_shapeOpti_pde} and \eqref{eq_defG_pde}, respectively, and define the homotopy $\mathcal H(\Omega,t)$ as in \eqref{eq_defH_aca}. This would mean that we just work with the averages of $u_G$ (solution to $e_G$) and $u_F$ (solution to $e_F$). In the approach defined here, we also average the equations, which is, in general, not the same. Along the same lines, one could also define the reduced cost function $\mathcal J_F(\Omega)$ of \eqref{eq_shapeOpti_pde} and do a convex combination with the unconstrained problem \eqref{eq_defG_aca}. However, we postulate that modeling a continuous transition between an artificial PDE constraint $e_G$ and the PDE of interest $e_F$ is beneficial for the arising homotopy path. We stress that a thorough comparison of different homotopies is beyond the scope of this article.
\end{remark}

\black
\subsection{A predictor-corrector scheme for shape optimization} \label{sec_predCorr_shape}
The procedure for solving shape optimization problems \eqref{eq_shapeOpti_aca} or \eqref{eq_shapeOpti_pde} using homotopies such as the ones defined in \eqref{eq_defHopt_aca} and \eqref{eq_defH_pde_reduced} is summarized in Algorithm \ref{algo_predCorr_shape_aca}. Starting out from a local minimizer $\Omega(0)$ of an auxiliary problem \eqref{eq_defG_aca} or \eqref{eq_defG_pde}, we proceed by finding stationary points of \eqref{eq_defHopt_aca} and \eqref{eq_defH_pde_reduced} for increasing values of the homotopy parameter $t$, i.e., by solving
\begin{align} \label{eq_dH_zero}
    d \mathcal H(\Omega(t), t)[\cdot] = 0
\end{align}
for $\Omega(t)$. We give details about the predictor step, the step size choice and the corrector steps in Sections \ref{sec_pred_shape}, \ref{sec_step_shape} and \ref{sec_corr_shape}, respectively.

\begin{algorithm}
    \begin{algorithmic}
    \STATE{Input: $\Omega^{(0)} \in \mathcal A, t^{(0)}=0, k=0$.}
    \STATE{Choose $\Delta t^{(0)}$}
    \WHILE{$t^{(k)}<1$}
        \STATE{Set $t^{(k+1)} = \mbox{min}(t^{(k)}+\Delta     t^{(k)}, 1).$ }
        \STATE{\textbf{Predictor:} $\tilde{\Omega}^{(k+1)} = \mathrm{Pred}^{(\npred)}(  \Omega^{(k)}, t^{(k)}, t^{(k+1)}-t^{(k)})$}
        \STATE{\textbf{Corrector:} Find stationary point $\Omega^{(k+1)}$ of $\mathcal H(\Omega, t^{(k+1)})$ starting from $\tilde{\Omega}^{(k+1)}$ by Algorithm \ref{algo_shapeNewton}.}
        \IF {corrector step successful}
            \STATE{Choose $\Delta t^{(k+1)}$}
            \STATE{$k \leftarrow k+1$}
        \ELSE
            \STATE{Reduce $\Delta t^{(k)}$}
        \ENDIF
    \ENDWHILE
    \end{algorithmic}
    \caption{Predictor-corrector for shape optimization}
    \label{algo_predCorr_shape_aca}
\end{algorithm}

\subsubsection{Predictors in shape optimization} \label{sec_pred_shape}

The results of Theorem \ref{thm_recursion_homo} carry over to this case by formally replacing differentiation with respect to $x$ by shape differentiation. Exemplarily, we give the equations defining the first, second and third order derivatives $\Omega'(t)$, $\Omega^{[2]}(t)$ and $\Omega^{[3]}(t)$ of the mapping $t\mapsto \Omega(t)$, cf. Corollary \ref{cor_predictors},
\begin{align} \label{eq_dtH_shape}
    d^2 \mathcal H[\Omega'(t), \cdot] =&  - d\mathcal H_t[\cdot], \\
    d^2 \mathcal H[\Omega^{[2]}(t), \cdot] =& \label{eq_dttH_shape}
    -  d^3 \mathcal H[\Omega'(t), \Omega'(t), \cdot] - 2d^2 \mathcal H_t[\Omega'(t), \cdot]  - d \mathcal H_{tt},\\
    \begin{split} \label{eq_dtttH_shape}
   d^2 \mathcal H[\Omega^{[3]}(t), \cdot] =& - d^4 \mathcal H[\Omega'(t), \Omega'(t), \Omega'(t), \cdot] - 3 d^3\mathcal H_t[\Omega'(t), \Omega'(t), \cdot] \\
    &- 3 d^3 \mathcal H[\Omega^{[2]}(t), \Omega'(t), \cdot]  - 3 d^2 \mathcal H_t[\Omega^{[2]}(t), \cdot] \\&- 3 d^2 \mathcal H_{tt}[\Omega'(t), \cdot] - d \mathcal H_{ttt}[\cdot],
    \end{split}
\end{align}
where we again omitted the arguments $(\Omega(t), t)$ of $\mathcal H$ for brevity. Again, note that the operator on the left hand side which is to be inverted is the same in all cases and thus, after discretization, the corresponding matrix has to be assembled (and, in the case of direct solvers, also factorized) only once for given $(\Omega(t^{(k)}), t^{(k)})$. Also note that this operator is the same as the one used in each iteration of Newton's method, cf. Section \ref{sec_corr_shape}. The computation of the higher order shape derivatives appearing in \eqref{eq_dtH_shape}--\eqref{eq_dtttH_shape} is discussed in Section \ref{sec_shapeDiff_higher}, see also \cite{GanglEtAlSAMO2020}. Note that systems \eqref{eq_dtH_shape}--\eqref{eq_dtttH_shape} have to be solved as discussed in Section \ref{sec_shapeNewton} and the predictors $\Omega'(t)$, $\Omega^{[2]}(t)$, $\Omega^{[3]}(t)$ are actually vector fields defined only on the boundary of $\Omega$ and have to be extended to the interior by \eqref{eq_extension}. Given $\npred$ derivatives of the domain with respect to the homotopy parameter $t$ and denoting the extension of the $i$-th derivative $\Omega^{[i]}(t)$ by $\hat \Omega^{[i]}(t)$, we can define the $\npred$-th order shape predictor at homotopy parameter $t^{(k)}$ by
\begin{align} \label{eq_pred_shape}
    \text{Pred}^{(\npred)}(\Omega^{(k)}, t^{(k)}, \Delta t^{(k)}) := \left( \text{id}  + \sum_{i=1}^\npred \frac{1}{i!}\left( \Delta t^{(k)} \right)^i \hat \Omega^{[i]}(t^{(k)}) \right) \left(\Omega^{(k)}\right).
\end{align}
\begin{remark} \label{rem_extension}
Note that, for $\npred>1$, it can be useful to, instead of extending each of the terms $\Omega^{[i]}$ individually, to apply the extension operator only to the weighted sum in \eqref{eq_pred_shape}. However, this entails that the extension has to be repeated every time when the homotopy step size $\Delta t^{(k)}$ is reduced.
\end{remark}

\subsubsection{Step size choice}  \label{sec_step_shape}
The step size choices presented in Sections \ref{sec_fixed_adap}--\ref{sec_step_agileadap} can be directly transfered to the setting of shape optimization. We mention that for the choice \eqref{eq_agile_step_size}, in our experiments the magnitude of the predictor of order $\npred+1$ is measured as $\| \Omega^{[\npred+1]} (t^{(k)}) \|_{L^2(\partial \Omega)^d}$.

\subsubsection{Shape-Newton method} \label{sec_corr_shape}
In the corrector step, given a homotopy value $t = t^{(k+1)}$ and starting out from a starting guess $\tilde \Omega^{(k+1)}$, we solve optimality condition \eqref{eq_dH_zero} by Newton's method as described in Algorithm \ref{algo_shapeNewton}.

In the unconstrained case, the Newton system is as given in \eqref{eq_shapeNewton_Vn} with $\mathcal J(\Omega)$ replaced by $\mathcal H(\Omega, t)$.

In the case of PDE-constrained shape optimization, the optimality condition \eqref{eq_dH_zero} can be written as
\begin{align}
    \begin{pmatrix}
    \partial_\Omega L_H(\Omega, u, p;t)[\cdot] \\
    \partial_u L_H(\Omega, u, p;t)[\cdot] \\
    \partial_p L_H(\Omega, u, 0;t)[\cdot]
    \end{pmatrix} = 0
\end{align}
for the Lagrangian $L_H(\Omega, u, p; t) := J_H(\Omega, u; t) + e_H(\Omega, u; t)(p)$. Assuming that the state $u$ and adjoint $p$ are approximated by finite element functions $u_h, p_h \in V_h$ with $V_h$ as defined in \eqref{eq_def_Vh}, i.e., $u_h(x) = \sum_{j=1}^N u_j \varphi_j(x)$, $p_h(x) = \sum_{j=1}^N p_j \varphi_j(x)$ with coefficient vectors $\underline u = (u_1,\dots, u_N)^\top, \underline p= (p_1,\dots, p_N)^\top \in \R^N$,  the corresponding discretized Newton system \eqref{eq_shapeNewton_Vn} reads
\begin{align}
\begin{pmatrix}
A_{\Omega,\Omega} & B_{\tau, \Omega} & A_{u, \Omega} & A_{p,\Omega} \\
B_{\tau,\Omega}^\top &  & & \\
A_{\Omega,u} & & A_{u, u} & A_{p,u} \\
A_{\Omega,p} & & A_{u, p} & \\
\end{pmatrix}
\begin{pmatrix}
    \underline V \\ \underline \xi \\ \underline u \\ \underline p
\end{pmatrix}
= -
\begin{pmatrix}
    \big[ \partial_\Omega L_H [\Phi_i] \big]_{i=1,\dots,d \Nbdy} \\
        \\
    \big[ \partial_u L_H [\varphi_l] \big]_{l=1,\dots,N} \\
    \big[ \partial_p L_H [\varphi_l] \big]_{l=1,\dots,N}
\end{pmatrix}
\end{align}
with the matrix entries
\begin{align*}
    \left(A_{\Omega, \Omega} \right)_{i,j} =& \partial_\Omega \partial_\Omega L_H [\Phi_i, \Phi_j], &
    \left(A_{u, \Omega} \right)_{i,k} =& \partial_u \partial_\Omega L_H [\Phi_i, \varphi_k],&
    \left(A_{p, \Omega} \right)_{i,k} =& \partial_p \partial_\Omega L_H [\Phi_i, \varphi_k], \\
    \left(A_{u, u} \right)_{l,k} =& \partial_u \partial_u L_H [\varphi_l, \varphi_k], &
    \left(A_{p, u} \right)_{l,k} =& \partial_p \partial_u L_H [\varphi_l, \varphi_k], &&
\end{align*}
for $i,j=1,\dots, d\Nbdy$ and $k,l=1,\dots,N$; moreover, $A_{\Omega,u}=A_{u,\Omega}^\top$, $A_{\Omega,p}=A_{p,\Omega}^\top$, $A_{u,p}=A_{p,u}^\top$ and $B_{\tau, \Omega}$ is as defined in \eqref{eq_def_BtauOmega}. Here, again, we skipped the arguments $(\Omega, u_h, p_h; t)$ of $L_H$.

\begin{remark} \label{rem_tolerance}
    For efficiency reasons, it seems reasonable to adapt the stopping criterion of the corrector step in the course of the homotopy. Since one is interested in accurate solutions of \eqref{eq_dH_zero} only for $t=1$, we suggest to adapt the stopping criterion by using $((1-t)r + t) \mathrm{tol}$ for a certain ratio $r\gg 1$ and the tolerance $\mathrm{tol}$ that is desired at $t=1$.
\end{remark}

\section{Homotopy methods for multi-objective shape optimization} \label{sec_pareto}
In multi-objective optimizaion one wants to simultaneously minimize multiple, possibly conflicting, cost functions. Since, in most cases, there does not exist a unique minimizer for all cost functions at the same time, one is interested in locally \textit{Pareto optimal} designs, i.e., designs which cannot be improved with respect to any of the given cost functions without deteriorating another cost function value. The set of all Pareto optimal points is usually called Pareto set and their image (under the application of the given cost functions) is called the Pareto front. Given the Pareto set, a designer can choose the most suitable solution for the particular problem at hand, often accounting for practical constraints such as manufacturability.

As a by-product of the chosen convex homotopy \eqref{eq_defH_aca}, the solutions to all intermediate problems $\mathcal H(\cdot, t) = 0$, $0 \leq t \leq 1$ with $\mathcal H$ defined by \eqref{eq_defH_aca}, correspond to locally Pareto optimal points of the bi-objective optimization problem to minimize $(\mathcal J_G(\Omega), \mathcal J_F(\Omega))^\top$ given by \eqref{eq_defG_aca}, \eqref{eq_shapeOpti_aca}. While this property seems to be of limited use since $\mathcal G$ represents an artificial, non-physical problem, we can, of course, also apply the methodology to connect two different physically relevant cost functions $\mathcal J_{F,1}(\Omega)$ and $\mathcal J_{F,2}(\Omega)$ by a homotopy $\mathcal H^{1,2}(\Omega, t) := (1-t)\mathcal J_{F,1}(\Omega) + t \mathcal J_{F,2}(\Omega)$. Starting out from a stationary point of $\mathcal J_{F,1}$ (which may itself have been obtained by a homotopy method) and $t=0$, the predictor-corrector method of Algorithm \ref{algo_predCorr_shape_aca} computes stationary points of $\mathcal H^{1,2}(\cdot, t)$ for different values of $t$, each of which corresponds to a point in the Pareto set.

We mention that the described concept of Pareto front tracing by predictor-corrector-type homotopy methods has previously been applied to multi-objective optimization \cite{Hillermeier2001JOTA, Hillermeier2001book, MartinSchuetze2018, schmidt2008pareto}. Notably, \cite{Hillermeier2001JOTA, Hillermeier2001book} treats the general case of an arbitrary number of cost functions, resulting in a generalized homotopy method with a multi-dimensional homotopy parameter. We also refer the interested reader to \cite{BoltenDoganayGottschalkKlamroth2021} for numerical integration techniques applied to the corresponding Davidenko differential equation \eqref{eq_Davidenko}. Beside the recent application \cite{CesaranoGanglSCEE} of some of the authors in the context nonlinear magnetostatics, to the best of our knowledge, the application in (free-form) shape optimization is novel.

\begin{remark}
 A common approach to obtaining a Pareto set by gradient-based optimization methods is the so-called weighted sum method where a locally Pareto optimal point is obtained as the result of a steepest descent method applied to a convex combination of given cost functions. Here, different points are obtained by varying the weight in the convex combination, i.e., a full optimization run is required in order to get one additional point on the Pareto front. With the homotopy method on the contrary, a neighboring point on the Pareto front can be obtained by solving just one corrector equation of the type \eqref{eq_shapeNewton_Vn} (with $\mathcal J(\Omega)$ replaced by $\mathcal H^{1,2}(\Omega, t)$), typically resulting in a significant speed-up of the procedure. Moreover, in contrast to the weighted sum method, the homotopy method is also able to determine the Pareto front if it is non-convex (more precisely: if the set in the objective space that is bounded by the Pareto front is non-convex) \cite{Hillermeier2001JOTA, Hillermeier2001book}.
\end{remark}

\begin{remark}
    The approach is also applicable in the case of PDE-constrained shape optimization. In this setting, however, the corresponding PDE constraints (if they are different) should not be mixed as in \eqref{eq_defH_pde}, but the convex combination of the corresponding reduced cost functions should be used.
\end{remark}

\begin{remark}
    We mention that, in practice, Pareto fronts may be non-connected and multiple starting points may be necessary to find different parts of the Pareto front.
\end{remark}

\section{Numerical results} \label{sec_numerics}
Finally, we present numerical results in which we illustrate different aspects of the method described above. We begin by illustrating the shape-Newton method proposed in Section \ref{sec_shapeNewton}, before showing results obtained by the homotopy method for shape optimization without and with PDE constraints. We conclude the section with a multi-objective shape optimization problem.

All numerical experiments are conducted within the finite element software NGSolve \cite{Schoeberl2014} using globally continuous, piecewise linear finite elements to represent both the deformations and also solutions of constraining PDEs. All arising systems of linear equations are solved by sparse direct solvers.
The corrector equations $\mathcal H(\cdot, t)=0$ are solved by Algorithm \ref{algo_predCorr_shape_aca} where $\text{tol}$ is chosen in dependence of the homotopy parameter, $\text{tol} = (1-t) 10^{-4} + t 10^{-10}$, since high accuracy is only required for larger values of $t$.
Concerning the discussion in Remark \ref{rem_extension}, we chose to compute an extension for the sum of all predictors rather than for each term individually.

The code used for obtaining the results presented below is available at \cite{shapeHomotopyZenodo2024}.

\subsection{Local convergence of shape-Newton algorithm} \label{sec_num_shapeNewton}
\begin{figure}
    \begin{tabular}{cc}
        \includegraphics[width=.45\textwidth]{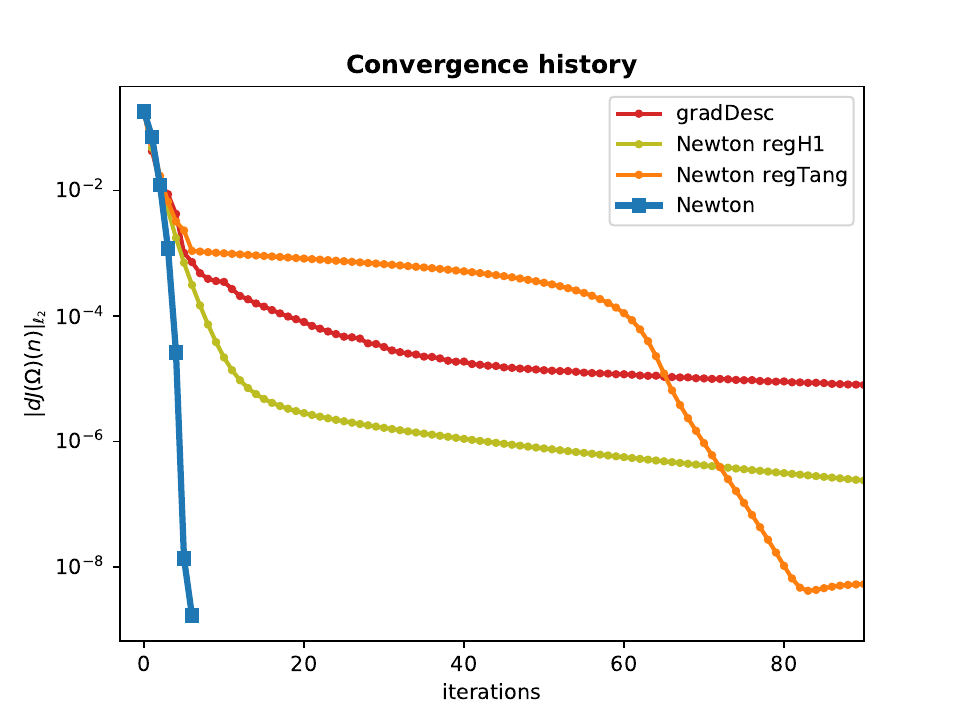} &
        \includegraphics[width=.45\textwidth]{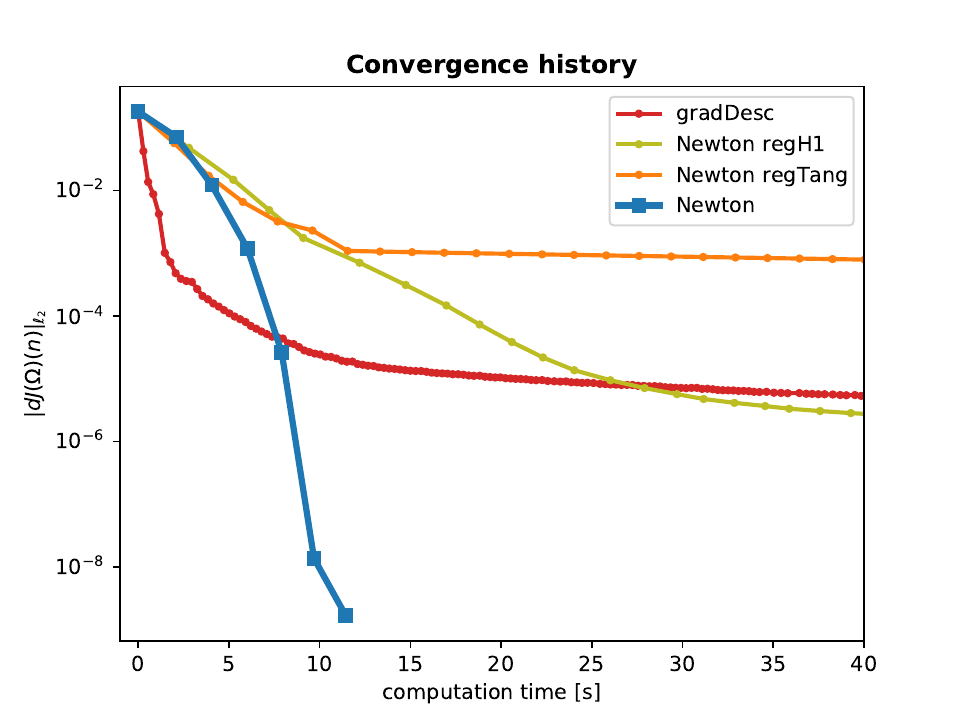} \\
        (a) & (b)
    \end{tabular}
    \caption{Convergence history of four different methods applied to the problem described in Section~\ref{sec_num_shapeNewton}: Shape derivative in normal direction $\| d\mathcal J(\Omega)(\vec{\mathbf n})\|_{\ell_2}$ as function of (a) iterations and (b) computation time.}
    \label{fig_ex8largerstepsizegradient}
\end{figure}

First, we numerically investigate the unregularized shape-Newton method proposed in Section \ref{sec_shapeNewton}, which will be used as a corrector in the experiments of the subsequent subsections. Since the method is only locally convergent, we apply it to a problem whose optimal solution is known and close to the initial design. More precisely, we consider the setting of Section \ref{sec_homo_shape_aca} with the integrand of \eqref{eq_shapeOpti_aca}, $f(x_1, x_2) = x_1^2/a^2 + x_2^2/b^2 - 1$, being a level set function of an ellipse centered at the origin with major and minor axes $a=1.25$, $b=1/a$, respectively. The initial design is the unit ball $\Omega^{(0)}=B_1(0)$ which is discretized by a triangular mesh with 3788 triangles and 1965 vertices and the integrand in \eqref{eq_defG_aca} is given by $\psi[\Omega^{(0)}](x_1, x_2) = x_1^2+x_2^2-1$. Figure \ref{fig_ex8largerstepsizegradient} shows the convergence history of the unregularized shape-Newton method in comparison with other first and second order methods applied to the same problem. In particular with compare with: (i) a regular shape gradient method where the descent direction is obtained by solving a problem like \eqref{eq_desc_direction} with a bilinear form of linear elasticity type \eqref{eq_extension}; (ii) a regularized shape-Newton approach where the shape Hessian in \eqref{eq_shapeNewton} is regularized by addition of an $H^1$ inner product, i.e. $\widetilde{d^2 \mathcal J}(\Omega)(\V,\W) = d^2 {\mathcal J}(\Omega)(\V,\W) +  \int_\Omega  \delta_{1,a}\partial \V : \partial \W + \delta_{1,b} \V \cdot \W \; \mbox{d}x$ as \cite{a_SC_2018a, SchmidtSchulz2023}; and (iii) the tangentially regularized shape Hessian $\widetilde{d^2 \mathcal J}(\Omega)(\V,\W) = d^2 {\mathcal J}(\Omega)(\V,\W) + \delta_2 \int_\Omega  (\V \cdot \tau) (\W \cdot \tau) \; \mbox{d}x$ proposed in \cite{GanglEtAlSAMO2020}. Of course the results for approaches (ii) and (iii) depend on the choices of the parameters $\delta_{1,a}$, $\delta_{1,b}$, $\delta_{2}$, which we chose empirically as $\delta_{1,a}=\delta_{1,b}=0.5$ and $\delta_2=250$ as these values yielded the best results based on multiple conducted experiments. Figures \ref{fig_ex8largerstepsizegradient}(a) and (b) show the evolution of the $\ell_2$-norm of the residual in normal direction $d \mathcal J(\Omega)(\vec{\mathbf n}):= [ d \mathcal J(\Omega)(\Phi_{2i} n_1(x_i) + \Phi_{2i+1} n_2(x_i))]_{i=1, \dots \Nbdy}$ in the course of the optimization iterations and as a function of computing time, respectively. Here $(n_1(x_i), n_2(x_i))^\top$ denotes the normal vector at the boundary point $x_i$ (obtained by averaging from the adjacent edges). Note that, in the discretized setting, the full shape derivative vector $[d \mathcal J(\Omega)(\Phi_i)]_{i=1,\dots,dN}$ may contain tangential components.
 It can be seen that the proposed unregularized Newton method can reach high accuracy within very few iterations. Moreover, for the given problem size, the method also performs very well in terms of computing time.

\begin{remark} \label{rem_poorgrad}
    The poor behavior of the shape gradient method can be attributed to the fact that the discretized shape derivative may contain spurious components corresponding to interior or tangential perturbations. In \cite{EtlingHerzogLoayzaWachsmuth2020}, the authors propose to allow only for mesh deformations that correspond to normal forces acting on the boundary of the shape which can significantly improve the reached accuracy and avoid mesh degeneracy. However, even in that setting, the number of necessary iterations remains high and the authors also propose a second order method based on this principle.
\end{remark}

\newcommand{\patheight}{output/example8normalnew/nPred0/}

\subsection{Application of homotopy method to simple shape optimization problem} \label{sec_num_ex9}
\newcommand{\pathnine}{output_bu_20240409/example9/}
\newcommand{\pathnineFirstTwo}{output_bu_20240409/example9/nPred2/stepSize-1.0/}
\newcommand{\pathnineFirstThree}{output_bu_20240409/example9/nPred3/stepSize-1.0/}
\newcommand{\pathnineFirstFour}{output_bu_20240409/example9/nPred4/stepSize-1.0/}
\newcommand{\pathnineFirstFive}{output_bu_20240409/example9/nPred5/stepSize-1.0/}
\newcommand{\pathnineSecondTwo}{output_bu_20240409/example9/nPred2/stepSize0.02/}
\newcommand{\pathnineSecondThree}{output_bu_20240409/example9/nPred3/stepSize0.02/}
\newcommand{\pathnineSecondFour}{output_bu_20240409/example9/nPred4/stepSize0.02/}
\newcommand{\pathnineSecondFive}{output_bu_20240409/example9/nPred5/stepSize0.02/}
\newcommand{\pathnineThirdTwo}{output_bu_20240409/example9/nPred2/stepSize100.02/}
\newcommand{\pathnineThirdThree}{output_bu_20240409/example9/nPred3/stepSize100.02/}
\newcommand{\pathnineThirdFour}{output_bu_20240409/example9/nPred4/stepSize100.02/}
\newcommand{\pathnineThirdFive}{output_bu_20240409/example9/nPred5/stepSize100.02/}
\begin{figure}
    \begin{tabular}{cc}
        \includegraphics[width=.4\textwidth, trim=150 200 150 200, clip]{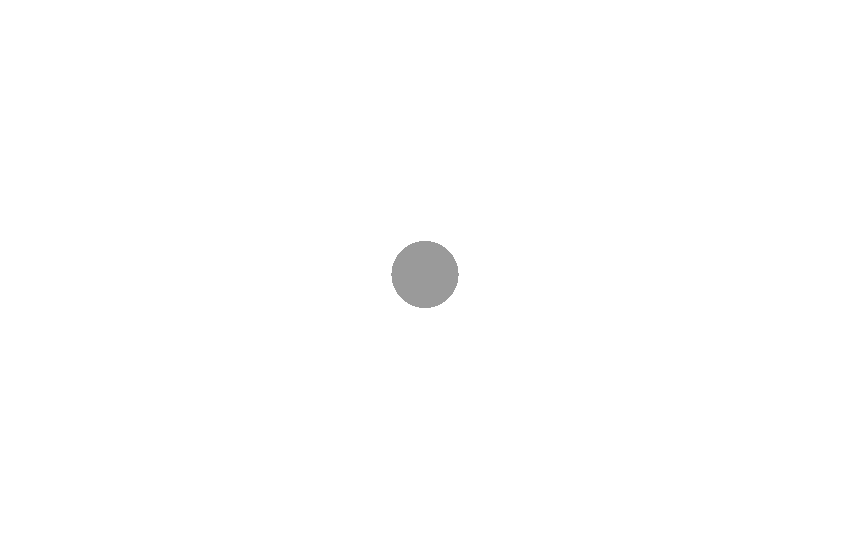} &
        \includegraphics[width=.4\textwidth, trim=150 200 150 200, clip]{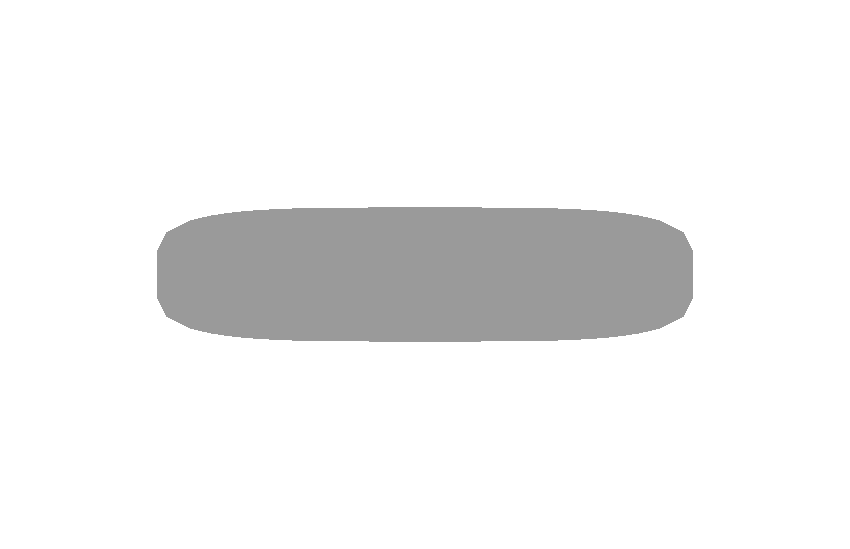} \\
        (a) & (b)
    \end{tabular}
    \caption{(a) Initial design and (b) optimal design for problem \eqref{eq_shapeOpti_aca} with function $f$ given by \eqref{eq_def_f_ex9}.}
    \label{fig_init_final_ex9}
\end{figure}
In this subsection, we illustrate the application of the homotopy methods discussed in Section \ref{sec_homo_shapeOpti} to the simple shape optimization problem \eqref{eq_shapeOpti_aca}. Also here, we consider as initial guess (an approximation of) the unit ball $\Omega^{(0)}=B_1(0)$ which is represented by the level set function $\psi[\Omega^{(0)}](x_1, x_2) = x_1^2+x_2^2-1$. For the function $f$ in \eqref{eq_shapeOpti_aca} we choose a function whose zero level set describes a large $p$-ellipse ,
\begin{align} \label{eq_def_f_ex9}
    f(x_1, x_2) = (x_1/a)^p + (x_2/b)^p - R^p,
\end{align}
with $p=4$, semi-axes $a=2$, $b=0.5$ and radius $R=4$, see Figure \ref{fig_init_final_ex9}. We discretized the initial design into a triangular mesh consisting of 2992 triangles and 1707 vertices. The example is chosen such that the optimal solution is relatively far away from the initial guess. As a result, Newton's method as described in Section \ref{sec_num_shapeNewton} fails to converge. Moreover, the evolution of a shape gradient method (with a line search for the step size) is very slow. After 10000 iterations and a runtime of 88 minutes, the residual was reduced only to about $5\cdot 10^{-3}$ (recall Remark \ref{rem_poorgrad}). % 5281.954631328583 %0.005428759109246723; 9169.48 seconds; see /home/pgangl/projectHannover/shapedwr/output/example9grad/nPred0/
For this example, we discuss the different potential ingredients to the predictor-corrector method described by Algorithm \ref{algo_predCorr_shape_aca}.

\begin{figure}
    \begin{center}
        \begin{tabular}{ccc}
            \hspace{-6mm}\includegraphics[width=.33 \textwidth]{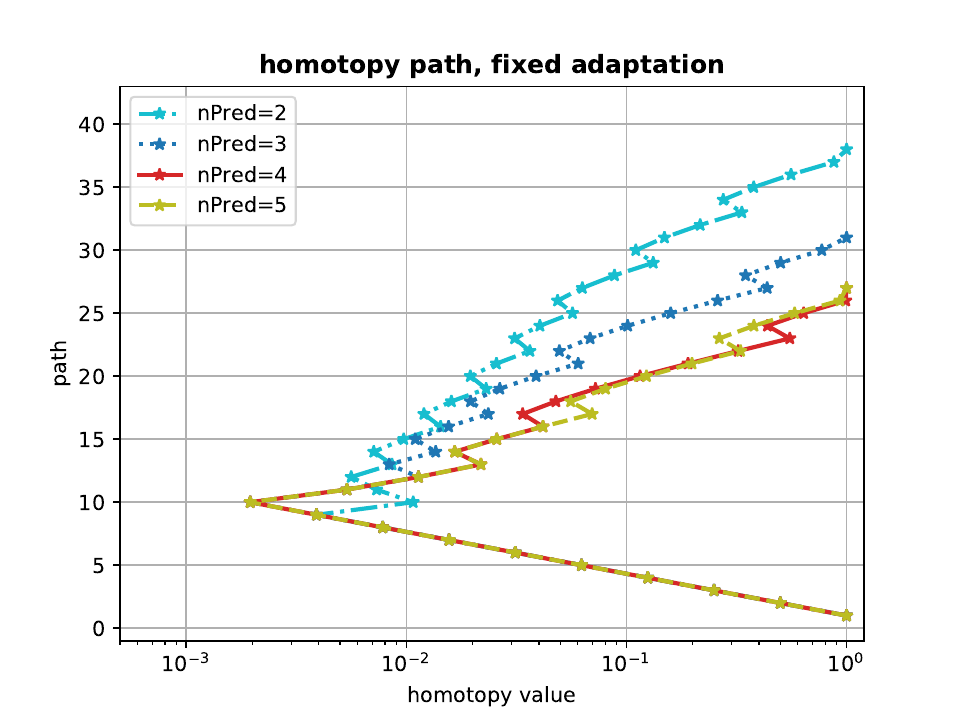}&
            \hspace{-3mm}\includegraphics[width=.33 \textwidth]{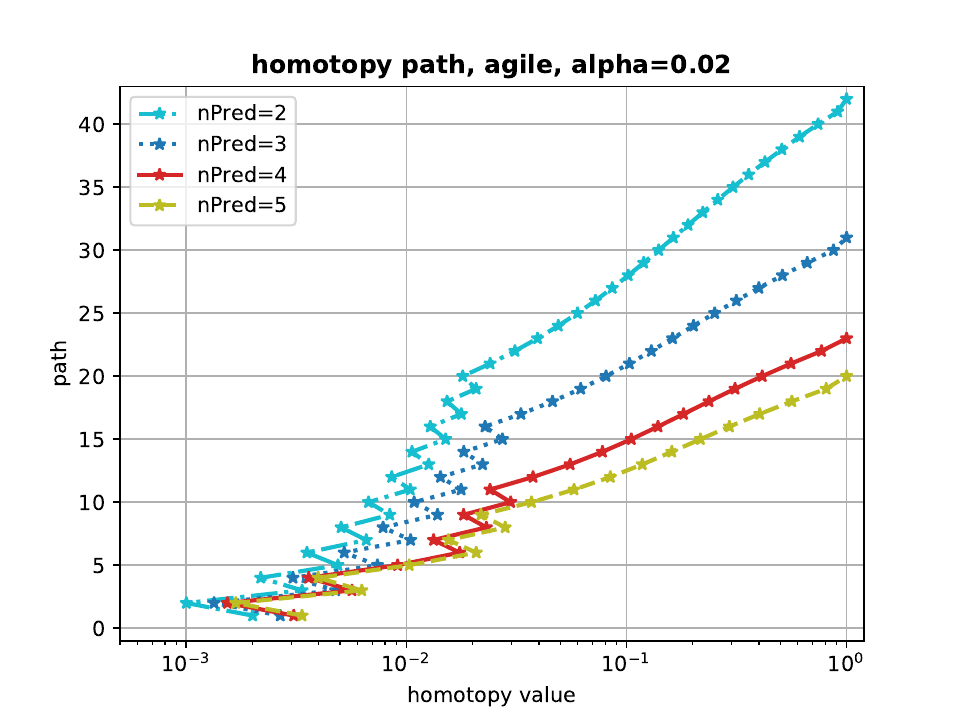}&
            \hspace{-3mm}\includegraphics[width=.33 \textwidth]{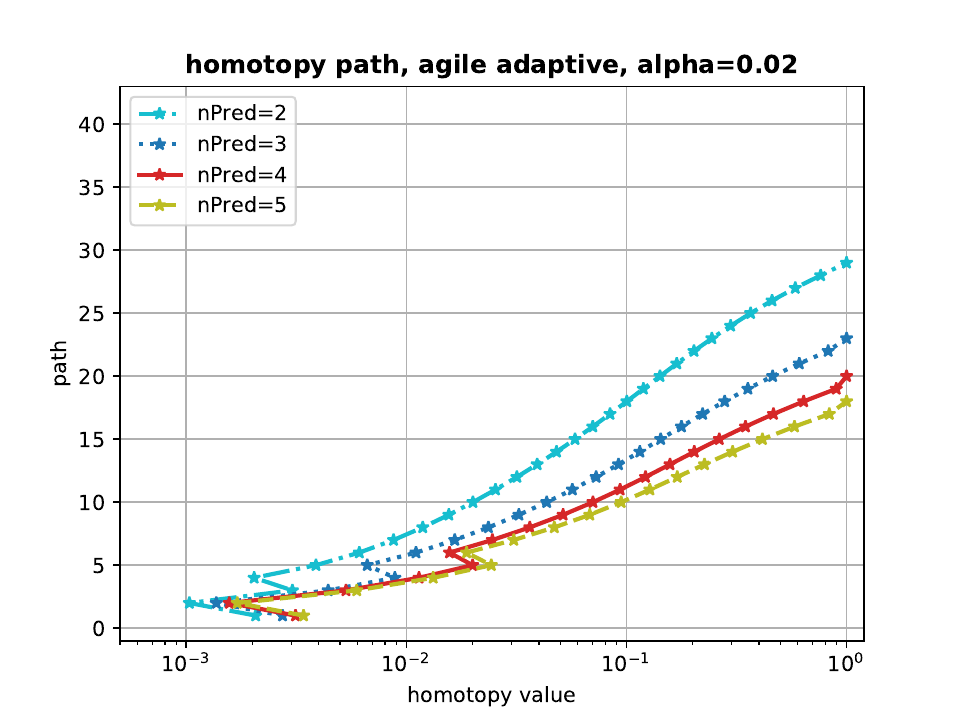}
        \end{tabular}
    \end{center}
    \caption{Homotopy paths for simple shape optimization example of Section \ref{sec_num_ex9} for $q=\mathrm{nPred}=2,3,4,5$. (a) Fixed step size adaptation with initial step size $\Delta t^{(0)}=1$ and update factors $\underline \gamma=0.5$, $\overline \gamma=1.75$. (b) Agile step size adaptation with $\alpha=0.02$. (c) Agile step size strategy with adaptation of $\alpha$ with initial choice $\alpha=0.02$ and update factors $\underline \alpha=0.5$, $\overline \alpha=1.1$.}
    \label{fig_ex9}
\end{figure}

\subsubsection{A first experiment with second order predictor}
In a first experiment, we use a second order predictor, $\npred=2$, together with the fixed adaptation step size choice described in Section \ref{sec_fixed_adap}. We start with $t^{(0)}=0$ and the initial step size $\Delta t^{(0)}=1$ and choose the reduction and increase factors $\underline \gamma = 0.5$ and $\overline \gamma = 1.75$, respectively. Note that $\Delta t^{(0)}=1$ is the most aggressive choice which would result in a pure Newton method in case of success.

After 39 visited homotopy values, out of which 22 were successful, a solution of $\mathcal H(\cdot,1)=0$ and thus a stationary point of \eqref{eq_shapeOpti_aca} is found, see Figure \ref{fig_init_final_ex9}(b).
% % \input{\pathnineFirstTwo t_all.txt}
% % \input{\pathnineFirstTwo t_succ.txt}
The solution required a total of
76
Newton steps and
42
solutions of predictor equations of the form \eqref{eq_dtH_shape} or \eqref{eq_dttH_shape}. While the computation time for all corrector steps was 323.21
seconds, the computation for the predictors was only 77.33
seconds. This difference can be explained by the fact that, for higher order predictors \eqref{eq_dttH_shape}, the system matrix does not have to be assembled and factorized again as it coincides with the system matrix of the previously computed first order predictor \eqref{eq_dtH_shape}.

\subsubsection{Effect of higher order predictor terms}
Next we investigate the effect of increasing the number of predictor terms in \eqref{eq_pred_shape}. We re-run the same experiment with $\npred=3, 4, 5$. Figure \ref{fig_ex9}(a) compares the different homotopy values attempted by our optimization runs. Depending on whether the corrector step for a given homotopy value was successful or not, the homotopy value is increased or decreased and the path proceeds to the right or left, respectively. Figure \ref{fig_ex9}(a) shows which homotopy values were passed for each of the four investigated optimization runs. The corresponding numbers of predictor and Newton steps as well as their timings are collected in the left column of Table \ref{tab_ex9}. It can be seen from Figure \ref{fig_ex9}(a) that, by employing higher order predictors, the number of necessary intermediate homotopy values can be decreased significantly. This way, the total number of Newton steps in the corrector stage gets smaller and smaller. On the other hand, however, the computational cost for computing these higher order predictors increases with their order. We emphasize that this is not primarily due to the additional solution of the corresponding predictor equations \eqref{eq_def_dtnx}, but much more due to the cost of assembling their right hand sides.

\begin{table}\footnotesize
    \caption{\noindent Section \ref{sec_num_ex9}: three different step size choices and predictors of order $\npred=2,3,4,5$, cf. Fig. \ref{fig_ex9}.
     \\
    ``Hsteps'' corresponds to visited homotopy values: $\mathrm{successful}$ $|$ $\mathrm{insuccessful}$ $|$ $\mathrm{\mathbf{total}}$;\\ ``Nsteps'' corresponds to number of solves of linear systems: predictor $|$ $\mathrm{corrector}$ $|$ $\mathrm{\mathbf{total}}$;\\ ``Ntime'' corresponds to computation time for linear systems in seconds: predictor $|$ $\mathrm{corrector}$ $|$ $\mathrm{\mathbf{total}}$
    }
    \begin{tabular}{cc|ccc}
    $\npred$&& fixed adaptation & agile $\alpha=0.02$ & agile $\alpha=0.02$ adaptive \\ \hline
    $2$ &Hsteps&
         22
$|$ 17
$|$ {\bf 39
} &
         33
$|$ 10
$|$ {\bf 43
} &
         28
$|$ 2
$|$ {\bf 30
} \\
    &Nsteps&
        {\it}$|$
        {}$|$
        {\bf 118
} &
        {\it 96
}$|$
        {80
}$|$
        {\bf 176
}  &
        {\it 81
}$|$
        {65
}$|$
        {\bf 146
} \\
    &Ntime &
        {\it}$|$
        {}$|$
        {\bf 400.54
} &
        {\it 140.75
}$|$
        {346.44
}$|$
        {\bf 487.19
} &
        {\it 118.24
}$|$
        {279.41
}$|$
        {\bf 397.65
}  \\ \hline
    $3$ &Hsteps&
         18
$|$ 14
$|$ {\bf 32
}    &
         24
$|$ 8
$|$ {\bf 32
} &
         22
$|$ 2
$|$ {\bf 24
}    \\
    &Nsteps&
        {\it 51
}$|$
        {59
}$|$
        {\bf 110
}  &
        {\it 92
}$|$
        {57
}$|$
        {\bf 149
}&
        {\it 84
}$|$
        {48
}$|$
        {\bf 132
}  \\
    &Ntime &
        {\it 73.48
}$|$
        {252.61
}$|$
        {\bf 326.09
} &
        {\it 141.17
}$|$
        {243.11
}$|$
        {\bf 384.28
} &
        {\it 128.18
}$|$
        {205.32
}$|$
        {\bf 333.51
}  \\ \hline
    $4$ &Hsteps&
        16
$|$ 12
$|$ {\bf 28
}    &
        19
$|$ 5
$|$ {\bf 24
}&
        19
$|$ 2
$|$ {\bf 21
}    \\
    &Nsteps&
        {\it 60
}$|$
        {46
}$|$
        {\bf 106
} &
        {\it 90
}$|$
        {42
}$|$
        {\bf 132
} &
        {\it 90
}$|$
        {42
}$|$
        {\bf 132
}  \\
    &Ntime &
        {\it 101.79
}$|$
        {215.76
}$|$
        {\bf 317.55
} &
        {\it 209.44
}$|$
        {212.32
}$|$
        {\bf 421.76
} &
        {\it 189.06
}$|$
        {185.89
}$|$
        {\bf 374.95
}  \\ \hline
    $5$ &Hsteps&
         16
$|$ 12
$|$ {\bf 28
}     &
         17
$|$ 4
$|$ {\bf 21
}  &
         17
$|$ 2
$|$ {\bf 19
  }    \\
    &Nsteps&
        {\it 75
}$|$
        {42
}$|$
        {\bf 117
}   &
        {\it 96
}$|$
        {37
}$|$
        {\bf 133
}  &
        {\it 96
}$|$
        {37
}$|$
        {\bf 133
}    \\
    &Ntime &
        {\it 178.74
}$|$
        {212.29
}$|$
        {\bf 391.03
}&
        {\it 285.08
}$|$
        {157.3
}$|$
        {\bf 442.38
}&
        {\it 284.7
}$|$
        {157.82
}$|$
        {\bf 442.51
} \\ \hline
    \end{tabular}
    \label{tab_ex9}
\end{table}

\subsubsection{Effect of step size choice}
\begin{figure}
    \centering
    \begin{tabular}{cc}
        \includegraphics[width=.4\textwidth]{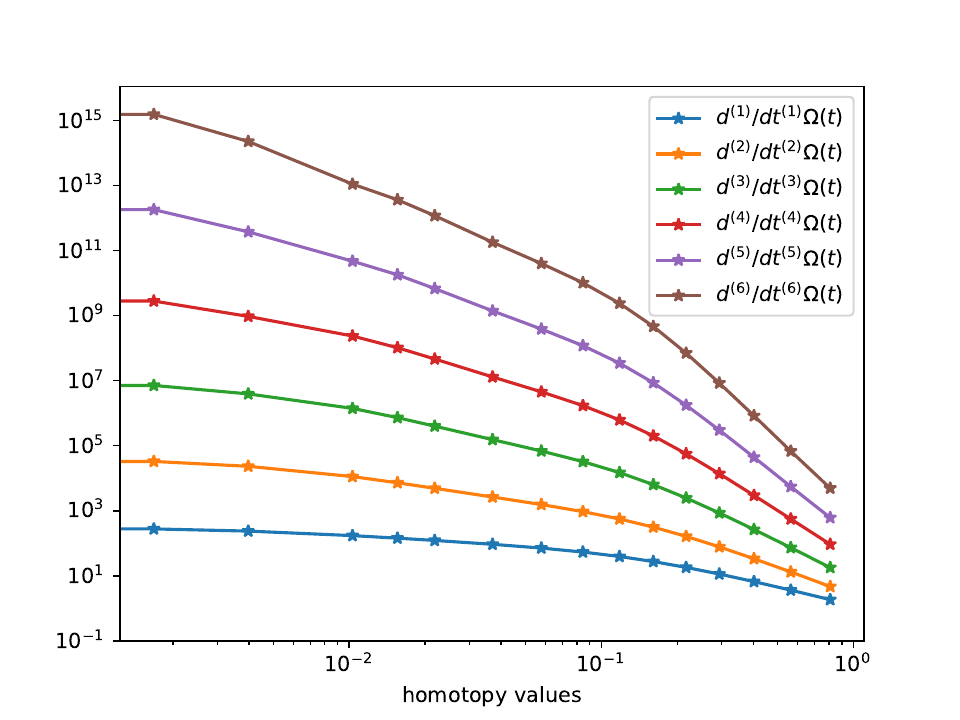} &
        \includegraphics[width=.4\textwidth]{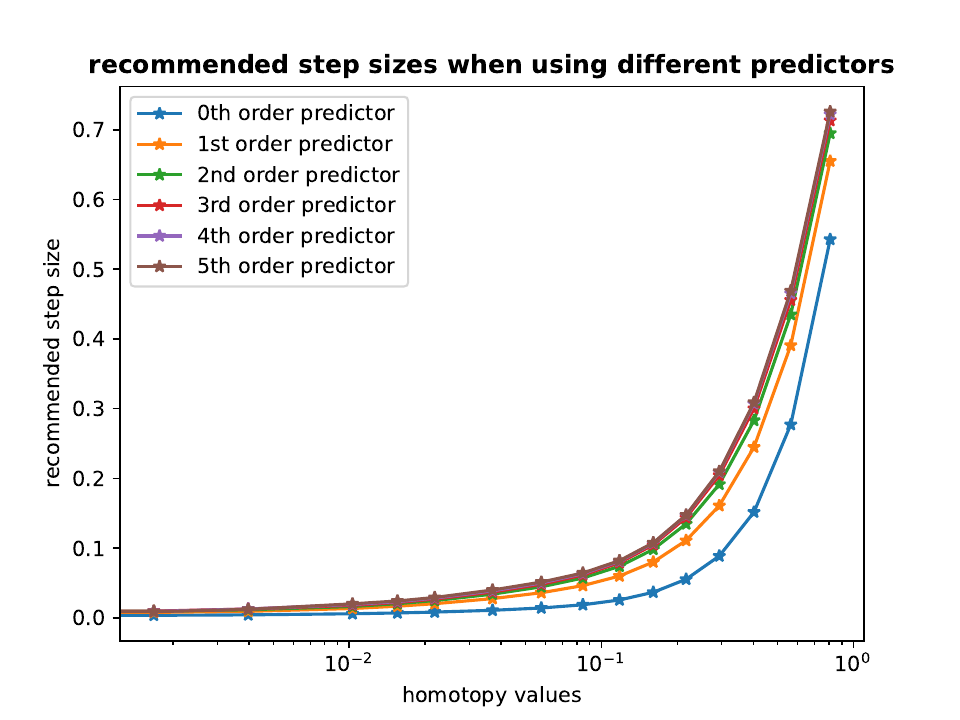} \\
        (a) & (b)
    \end{tabular}
    \caption{Example of Section \ref{sec_num_ex9}: (a) $L^2(\partial \Omega)$-norms of derivatives of $t \mapsto \Omega(t)$ of different orders (see, e.g., \eqref{eq_dtH_shape}--\eqref{eq_dtttH_shape}). (b) Recommended step sizes according to \eqref{eq_agile_step_size} for $\alpha=0.02$.}
    \label{fig_curvatures_ex9}
\end{figure}

Let us now have a closer look at the step size choices discussed in Sections \ref{sec_homo_stepsize} and \ref{sec_step_shape}. Note that, for the fixed adaptation strategy according to Section \ref{sec_fixed_adap} used in the experiments above, one has to choose three parameters (initial step size $\Delta t^{(0)}$ as well as the increase and decrease parameters $\overline \gamma, \underline \gamma$). The agile step size choice, on the other hand, automatically chooses the step size according to the Taylor remainder of the homotopy curve given a predictor of certain order. Thus, this method requires to set only one parameter $\alpha$ which is related to the allowed deviation from the homotopy curve. Application of this agile step size rule with $\alpha=0.02$ is illustrated in Figure \ref{fig_ex9}(b). When comparing with the result in Figure \ref{fig_ex9}(a), it can be seen that the method can automatically detect the order of the initial step size $\Delta t^{(0)}$. For instance, for $\npred=2$, the proposed step size has to be halvened only once compared to eight times in Figure \ref{fig_ex9}(a) when starting from the value $\Delta t^{(0)}=1$. The problem at hand requires comparably small step sizes in the beginning while larger step sizes can be made later on. Figure \ref{fig_curvatures_ex9}(a) shows the $L^2(\partial \Omega)$-norms of the different predictors $\Omega^{[j]}(t_k)$ for the derivative orders $j=1, \dots, 6$ for the homotopy values $t_k$ arising in the experiment of Figure \ref{fig_ex9}(b). Figure \ref{fig_curvatures_ex9}(b) shows the step sizes to be taken according to this rule as a result of \eqref{eq_agile_step_size}. Note that the step size rule automatically detects that small steps should be taken in the initial phase and larger steps are allowed later in the solution process.

While this agile step size choice has the advantage of adapting the step size to the properties of the physical problem, it still depends on the parameter $\alpha$. As it can be seen from Figure \ref{fig_ex9}(b), if the value of $\alpha$ is chosen too conservatively, the procedure may take many iterations. Therefore, we also investigate the strategy proposed in Section \ref{sec_step_agileadap} and increase or decrease the value of $\alpha$ by constant factors $\overline \alpha>1>\underline \alpha$ after every successful or insuccessful homotopy step, respectively. As it can be seen from Figure \ref{fig_ex9}(c) where we used $\overline \alpha=1.1$, $\underline \alpha=0.5$, the total number of homotopy steps can be reduced in this way. We remark, however, that for the two strategies followed in Figures \ref{fig_ex9}(b) and (c), an additional predictor term has to be computed. Thus, their applicability depends on whether an additional predictor term can be obtained without too much additional computational effort, see also the discussion in Remark \ref{rem_discuss_cost_predictores}. A comparison of the three strategies in terms of number of predictor or Newton solves and computation time can be found in Table \ref{tab_ex9}.

\begin{remark} \label{rem_meshqual}
    When shape optimization is performed based on mesh deformations, the aspect of mesh quality plays an important role, in particular in the presence of large deformations.
    Since the computation of mesh quality preserving descent directions is an active field of research on its own, we largely exclude this aspect from the discussion here and just mention that the mesh of Figure \ref{fig_ex9}(b) is distorted, but not degenerated. In order to deal with large deformations, one could employ local mesh refinement or remeshing techniques~\cite{AllaireDapogny2021}.
\end{remark}

\subsection{Application of homotopy method to PDE-constrained shape optimization} \label{sec_num_ex105}
\begin{figure}\centering
    \begin{tabular}{cc}
        \includegraphics[width=.4\textwidth, trim=100 70 100 70, clip]{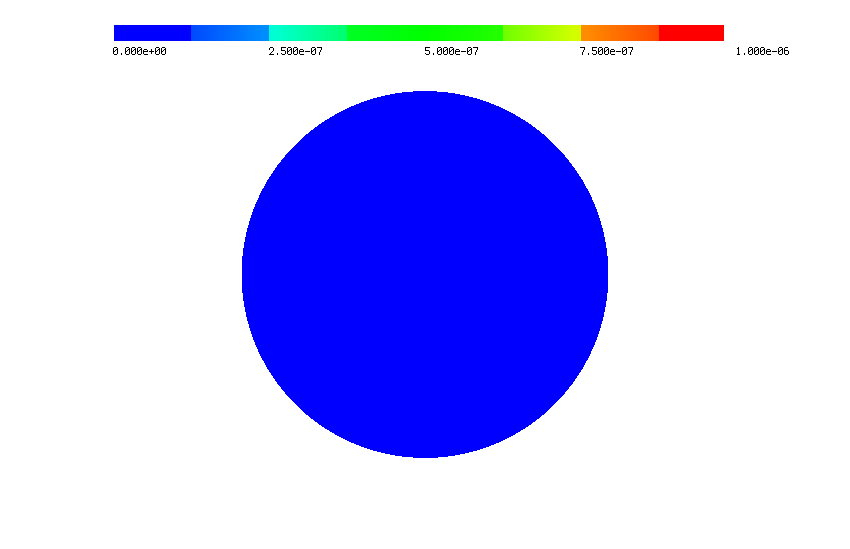} &
        \includegraphics[width=.4\textwidth, trim=100 70 100 70, clip]{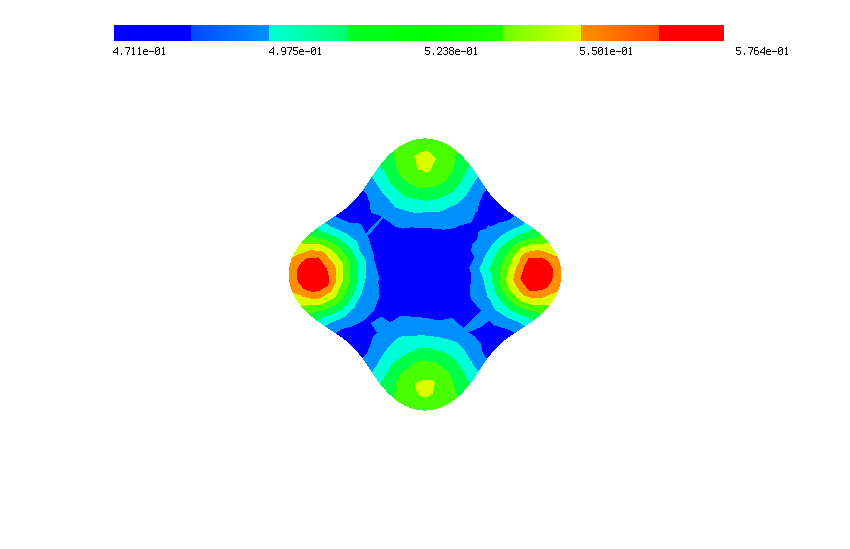} \\
        (a) & (b) \\
    \end{tabular}
    \caption{(a) Initial design and (b) final design for problem \eqref{eq_shapeOpti_pde} with $F$ and $e_F(\Omega, u)$ given by \eqref{eq_def_F_ex105} and \eqref{eq_def_eF_ex105}, respectively.}
    \label{fig_init_final_ex105}
\end{figure}

\begin{figure}
    \begin{center}
        \begin{tabular}{ccc}
            \hspace{-6mm}\includegraphics[width=.33 \textwidth]{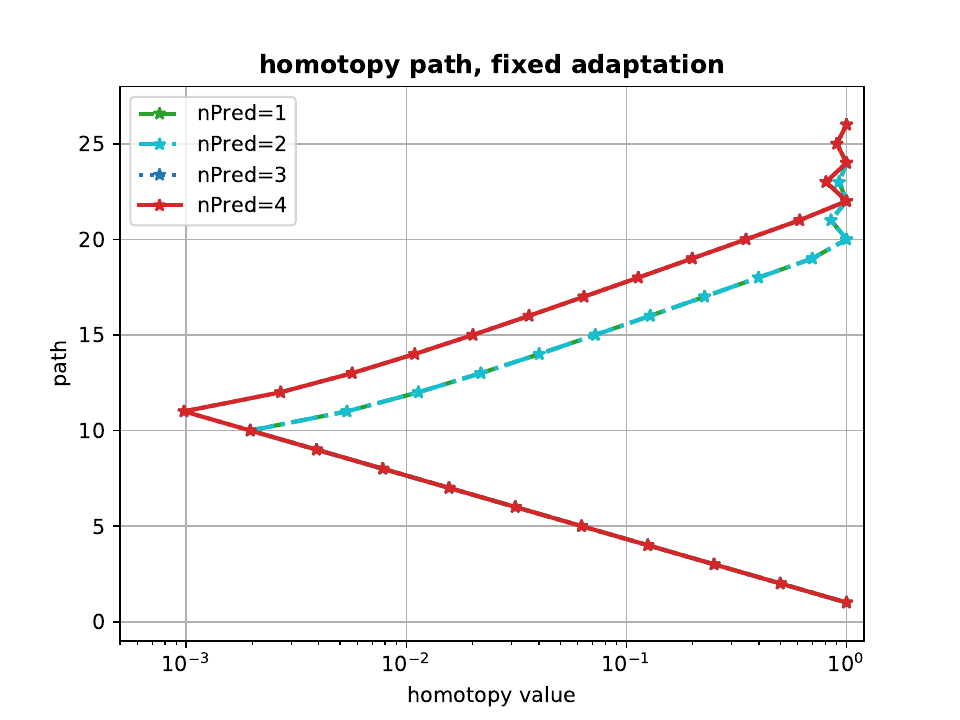}&
            \hspace{-3mm}\includegraphics[width=.33 \textwidth]{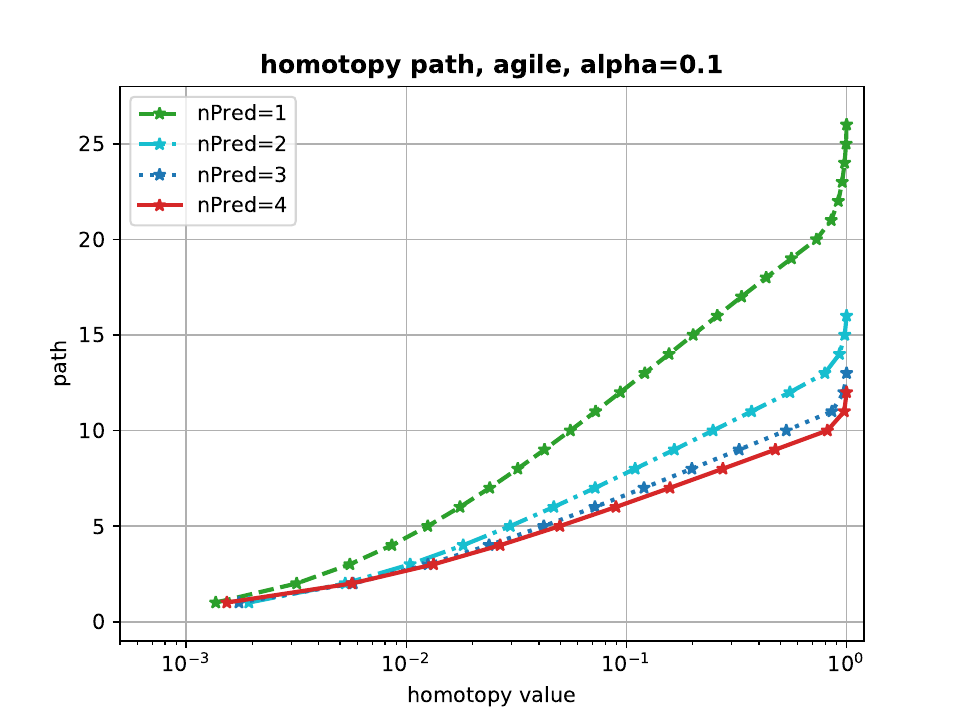}&
            \hspace{-3mm}\includegraphics[width=.33 \textwidth]{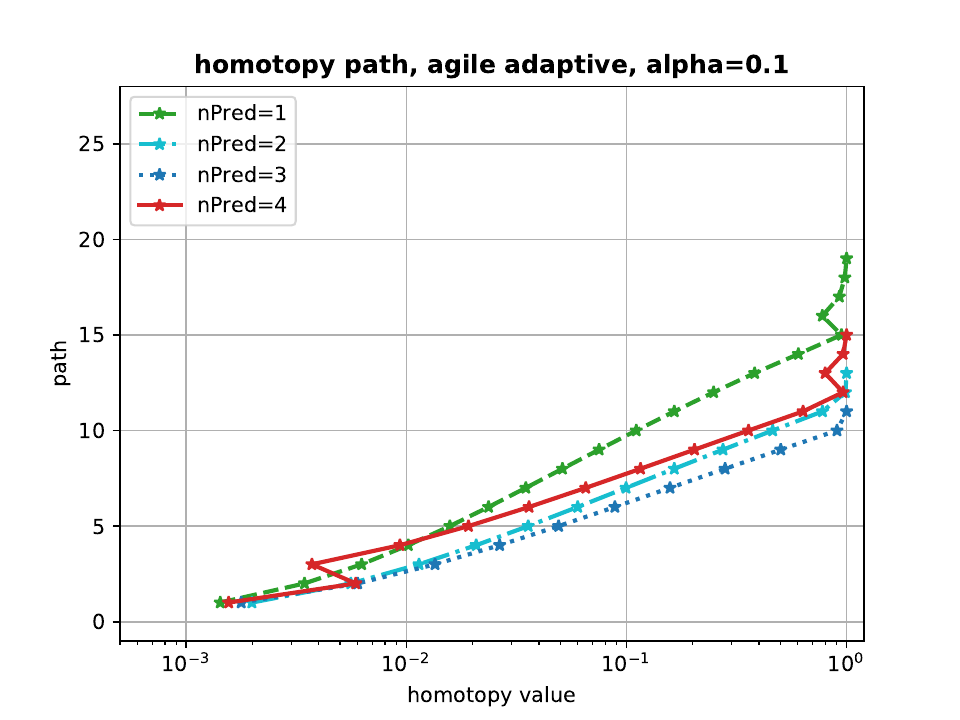} \\
            (a) & (b) & (c)
        \end{tabular}
    \end{center}
    \caption{Homotopy paths for PDE-constrained example of Section \ref{sec_num_ex105} for $q=\mathrm{nPred}=1,2,3,4$.  (a) Fixed step size adaptation with initial step size $\Delta t^{(0)}=1$ and update factors $\underline \gamma=0.5$, $\overline \gamma=1.75$. (b) Agile step size adaptation with $\alpha=0.1$. (c) Agile step size strategy with adaptation of $\alpha$ with initial choice $\alpha=0.1$ and update factors $\underline \alpha=0.5$, $\overline \alpha=1.1$.}
    \label{fig_ex105}
\end{figure}

\newcommand{\pathoneOfive}{output_bu_20240409/example105/}
Next, we consider the same strategies for a PDE-constrained shape optimization problem of the form \eqref{eq_shapeOpti_pde} with $V=H^1(\Omega)$,
\begin{align} \label{eq_def_F_ex105}
    F(\Omega, u) = \int_\Omega u(x) \; \mbox{d}x,
\end{align}
and the semilinear PDE constraint
\begin{align} \label{eq_def_eF_ex105}
    \langle e_F(\Omega, u), v \rangle = \int_\Omega \left[ \lambda(x) \nabla u \cdot \nabla v + u^3 v -f_\text{clov}v \right]\; \mbox{d}x, \; v \in H^1(\Omega),
\end{align}
with $\lambda(x_1, x_2) = 1/(1+x_1^2)$ and the right-hand side function whose zero level set resembles a clover shape
\begin{align} \label{eq_def_f_ex105}
\begin{aligned}
f_\text{clov}(x_1, x_2) =& \left(\sqrt{ (x_1 - a)^2 + b  x_2^2} - 1\right)
            \left(\sqrt{(x_1 + a)^2 + b  x_2^2} - 1\right) \\
            &\left(\sqrt{b x_1^2 + (x_2 - a)^2} - 1\right)
            \left(\sqrt{b x_1^2 + (x_2 + a)^2} - 1\right) - \epsilon,
\end{aligned}
\end{align}
where $a=0.8$, $b=2$ and $\epsilon=0.01$. We choose the homotopy defined by \eqref{eq_defH_pde} with the simple problem \eqref{eq_defG_pde_concrete} and the initial shape $\Omega^{(0)}=B_{2.5}(0)$ being the ball of radius $2.5$ centered at the origin with the corresponding initial level set function $\psi[\Omega^{(0)}](x_1,x_2) = x_1^2+x_2^2-2.5^2$. We used a mesh consisting of 650 triangles and 384 vertices. The initial design and the final design with the solution of the state equation are depicted in Figure \ref{fig_init_final_ex105}. The homotopy paths for the same three different step size choices as discussed in Section \ref{sec_num_ex9} are shown in Figure \ref{fig_ex105} and corresponding numbers of predictor and Newton solves as well as timings are given in Table \ref{tab_ex105}.

\begin{remark} \label{rem_discuss_cost_predictores}
    As it can be seen from Table \ref{tab_ex105}, the computation time for computing predictors grows fast with their order. This is mostly due to the assembly of the right hand sides of \eqref{eq_def_dtnx} where automated differentiation is used and could be reduced if higher order derivatives were provided in closed form. While the focus of this study is the investigation of the effects of higher order predictors, we mention that efficiency could be improved, e.g., by constructing predictors at $t^{(k)}$ by extrapolation from previously visited homotopy values $(t^{(j)}, x(t^{(j)}))$, $j \leq k$ as in multi-step methods for ordinary differential equations.
\end{remark}

\newcommand{\pathoneOfiveFirstTwo}{output_bu_20240409/example105/nPred1/stepSize-1.0/}
\newcommand{\pathoneOfiveFirstThree}{output_bu_20240409/example105/nPred2/stepSize-1.0/}
\newcommand{\pathoneOfiveFirstFour}{output_bu_20240409/example105/nPred3/stepSize-1.0/}
\newcommand{\pathoneOfiveFirstFive}{output_bu_20240409/example105/nPred4/stepSize-1.0/}
\newcommand{\pathoneOfiveSecondTwo}{output_bu_20240409/example105/nPred1/stepSize0.1/}
\newcommand{\pathoneOfiveSecondThree}{output_bu_20240409/example105/nPred2/stepSize0.1/}
\newcommand{\pathoneOfiveSecondFour}{output_bu_20240409/example105/nPred3/stepSize0.1/}
\newcommand{\pathoneOfiveSecondFive}{output_bu_20240409/example105/nPred4/stepSize0.1/}
\newcommand{\pathoneOfiveThirdTwo}{output_bu_20240409/example105/nPred1/stepSize100.1/}
\newcommand{\pathoneOfiveThirdThree}{output_bu_20240409/example105/nPred2/stepSize100.1/}
\newcommand{\pathoneOfiveThirdFour}{output_bu_20240409/example105/nPred3/stepSize100.1/}
\newcommand{\pathoneOfiveThirdFive}{output_bu_20240409/example105/nPred4/stepSize100.1/}

\begin{table}\footnotesize
    \caption{\noindent Section \ref{sec_num_ex105}: three different step size choices and predictors of order $\npred=1,2,3,4$, cf. Fig. \ref{fig_ex105}.
     \\
    ``Hsteps'' corresponds to visited homotopy values: $\mathrm{successful}$ $|$ $\mathrm{insuccessful}$ $|$ $\mathrm{\mathbf{total}}$;\\ ``Nsteps'' corresponds to number of solves of linear systems: predictor $|$ $\mathrm{corrector}$ $|$ $\mathrm{\mathbf{total}}$;\\ ``Ntime'' corresponds to computation time for linear systems in seconds: predictor $|$ $\mathrm{corrector}$ $|$ $\mathrm{\mathbf{total}}$
}
    \begin{tabular}{cc|ccc}
    $\npred$&& fixed adaptation & agile $\alpha=0.1$ & agile $\alpha=0.1$ adaptive \\ \hline
    $1$ &Hsteps&
         14
$|$ 11
$|$ {\bf 25
} &
         27
$|$ 0
$|$ {\bf 27
} &
         19
$|$ 1
$|$ {\bf 20
} \\
    &Nsteps&
        {\it 13
}$|$
        {61
}$|$
        {\bf 74
} &
        {\it 52
}$|$
        {70
}$|$
        {\bf 122
}  &
        {\it 36
}$|$
        {55
}$|$
        {\bf 91
} \\
    &Ntime &
        {\it 18.43
}$|$
        {84.43
}$|$
        {\bf 102.86
} &
        {\it 43.46
}$|$
        {88.75
}$|$
        {\bf 132.2
} &
        {\it 30.55
}$|$
        {71.37
}$|$
        {\bf 101.92
}  \\ \hline
    $2$ &Hsteps&
         14
$|$ 11
$|$ {\bf 25
}    &
         17
$|$ 0
$|$ {\bf 17
} &
         14
$|$ 0
$|$ {\bf 14
}    \\
    &Nsteps&
        {\it 26
}$|$
        {48
}$|$
        {\bf 74
}  &
        {\it 48
}$|$
        {43
}$|$
        {\bf 91
}&
        {\it 39
}$|$
        {38
}$|$
        {\bf 77
}  \\
    &Ntime &
        {\it 22.52
}$|$
        {60.43
}$|$
        {\bf 82.95
} &
        {\it 86.15
}$|$
        {55.66
}$|$
        {\bf 141.81
} &
        {\it 81.21
}$|$
        {60.91
}$|$
        {\bf 142.13
}  \\ \hline
    $3$ &Hsteps&
        15
$|$ 12
$|$ {\bf 27
}    &
        14
$|$ 0
$|$ {\bf 14
}&
        12
$|$ 0
$|$ {\bf 12
}    \\
    &Nsteps&
        {\it 42
}$|$
        {45
}$|$
        {\bf 87
} &
        {\it 52
}$|$
        {35
}$|$
        {\bf 87
} &
        {\it 44
}$|$
        {33
}$|$
        {\bf 77
}  \\
    &Ntime &
        {\it 81.18
}$|$
        {60.52
}$|$
        {\bf 141.7
} &
        {\it 650.94
}$|$
        {45.33
}$|$
        {\bf 696.27
} &
        {\it 560.82
}$|$
        {48.25
}$|$
        {\bf 609.08
}  \\ \hline
    $4$ &Hsteps&
         15
$|$ 12
$|$ {\bf 27
}     &
         13
$|$ 0
$|$ {\bf 13
}  &
         14
$|$ 2
$|$ {\bf 16
  }    \\
    &Nsteps&
        {\it 56
}$|$
        {41
}$|$
        {\bf 97
}   &
        {\it 60
}$|$
        {34
}$|$
        {\bf 94
}  &
        {\it 65
}$|$
        {39
}$|$
        {\bf 104
}    \\
    &Ntime &
        {\it 705.69
}$|$
        {60.44
}$|$
        {\bf 766.13
}&
        {\it 18587.18
}$|$
        {43.02
}$|$
        {\bf 18630.2
}&
        {\it 20994.39
}$|$
        {47.63
}$|$
        {\bf 21042.02
} \\ \hline
    \end{tabular}
    \label{tab_ex105}
\end{table}

\subsection{Pareto front tracing in shape optimization} \label{sec_num_pareto}
Finally, we illustrate the application of the proposed homotopy method for tracing Pareto optimal shapes. For that purpose, we consider three instances of the simple academic shape optimization problem \eqref{eq_shapeOpti_aca} with integrands
\begin{align*}
    f_1(x_1, x_2) = \frac{x_1^2}{a^2} + \frac{x_2^2}{b^2} - 2^2, &&
    f_2(x_1, x_2) = f_{\text{clov}}(x_1, x_2), &&
    f_3(x_1, x_2) = \frac{x_1^2}{b^2} + \frac{x_2^2}{a^2} - 2^2,
\end{align*}
with $a=1.3$, $b=1/a$ and $f_{\text{clov}}$ as given in \eqref{eq_def_f_ex105}. We define the convex combination of the three corresponding cost functions $\mathcal J_{F,i}(\Omega) = \int_\Omega f_i(x) \, \mbox{d}x$, $i=1,2,3$, as $\mathcal J_\text{conv}[s_1, s_2](\Omega) = s_1 \mathcal J_{F,1}(\Omega) + s_2 \mathcal J_{F,2}(\Omega) + (1-s_1-s_2) \mathcal J_{F,3}(\Omega)$, for $s_1, s_2 \geq0, s_1+s_2 \leq 1$. In order to trace along the Pareto front, which is a two-dimensional surface in $\R^3$, we define the three homotopies
\begin{align*}
    \tilde{ \mathcal H}^{1,2}(\Omega, t, \delta):=& (1-t) \mathcal J_{\text{conv}}[1-2 \delta, \delta](\Omega) + t \mathcal J_{\text{conv}}[\delta, 1-2 \delta](\Omega) \\
    \tilde{ \mathcal H}^{2,3}(\Omega, t, \delta):=& (1-t) \mathcal J_{\text{conv}}[\delta, 1-2 \delta](\Omega) + t \mathcal J_{\text{conv}}[\delta, \delta](\Omega) \\
    \tilde{ \mathcal H}^{3,1}(\Omega, t, \delta):=& (1-t) \mathcal J_{\text{conv}}[\delta, \delta](\Omega) + t \mathcal J_{\text{conv}}[1-2 \delta, \delta](\Omega)
\end{align*}
which are parametrized by $t \in [0,1]$ and $\delta \in [0,\frac13]$. Note that, for $\delta=0$, the homotopy $\tilde{\mathcal H}^{i,j}$ is just the convex homotopy between $\mathcal J_{F,i}$ and $\mathcal J_{F,j}$.

In our numerical experiments, we chose $\delta \in \{0, 0.1, 0.2, 0.3\}$ and for each fixed $\delta$ followed the homotopies $\tilde{\mathcal H}^{1,2}$, $\tilde{\mathcal H}^{2,3}$ and $\tilde{\mathcal H}^{3,1}$ from $t=0$ to $t=1$. Since, here, all intermediate solutions of $\tilde{\mathcal H}^{i,j}(\Omega, t, \delta)=0$ are of interest as they represent Pareto optimal shapes, we always use the full desired tolerance in the corrector stage, i.e., Remark \ref{rem_tolerance} does not apply here.
In order to solve $\tilde{\mathcal H}^{1,2}(\Omega,0, \delta)=0$ for any given $\delta$, we used the homotopy defined by \eqref{eq_defH_aca} with the simple problem \eqref{eq_defG_aca} related to the given initial design. The solutions of $\tilde{\mathcal H}^{2,3}(\Omega,0, \delta)=0$ and $\tilde{\mathcal H}^{3,1}(\Omega,0, \delta)=0$ are trivial starting out from solutions of $\tilde{\mathcal H}^{1,2}(\Omega,1, \delta)=0$ and $\tilde{\mathcal H}^{2,3}(\Omega,1, \delta)=0$, respectively. We used second order predictors and the agile step size strategy of Section \ref{sec_step_agile} with $\alpha=0.1$. In addition, in order to promote an equi-spaced Pareto front, we imposed a maximum distance in objective space of any two subsequent points in a homotopy run $\tilde{\mathcal H}^{i,j}$, i.e., we rejected points even if the corrector was successful and continued with a reduced step size.

\begin{figure}
    \begin{center}
        \includegraphics[width=.75 \textwidth, trim=100 150 50 200, clip]{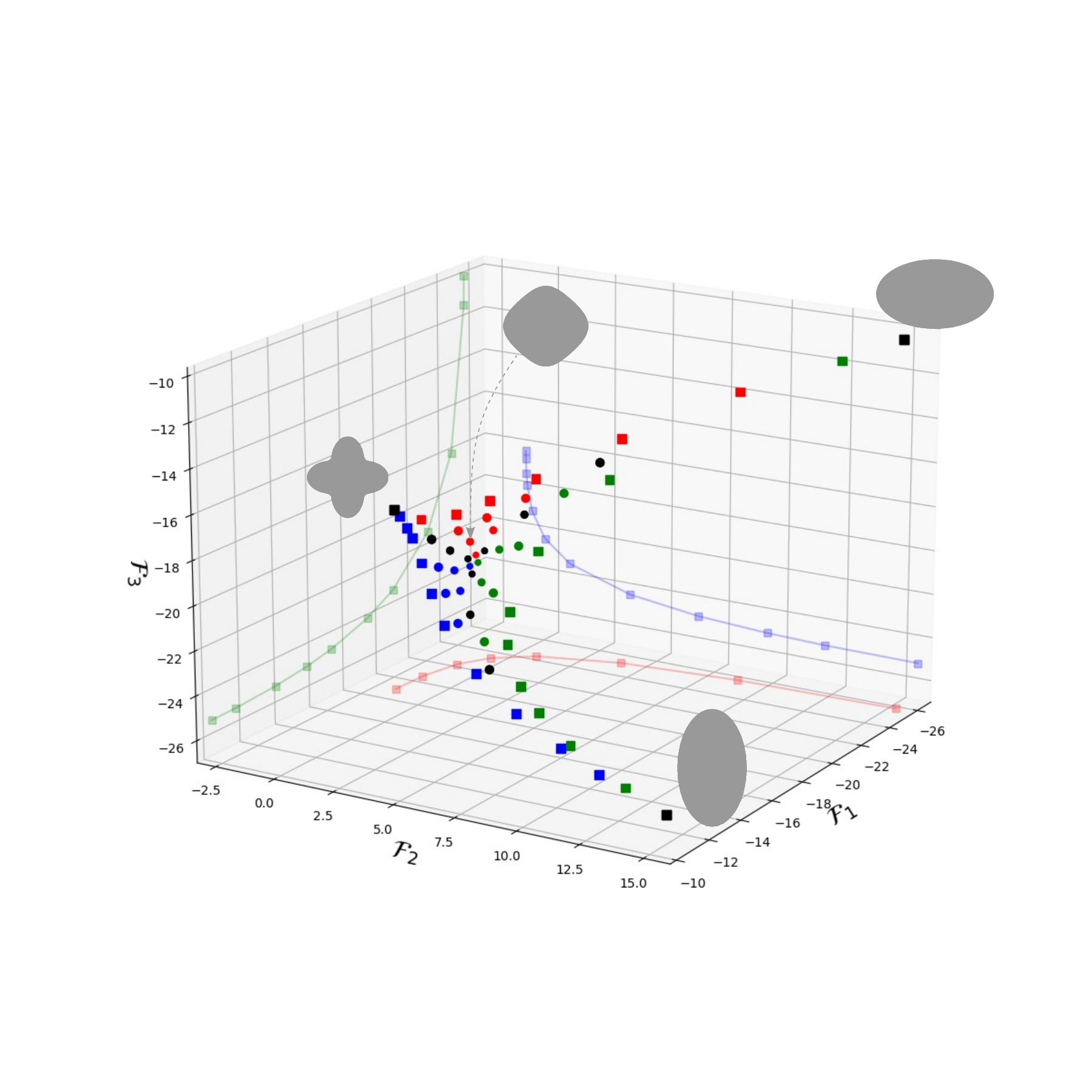}
    \end{center}
    \caption{Pareto front for three-objective shape optimization problem obtained by Pareto method of Section \ref{sec_num_pareto}. Black squares represent solutions to single-objective optimization problems; red, blue and green squares represent solutions to bi-objective optimization problems.}
    \label{fig_pareto}
\end{figure}

The resulting Pareto front is depicted in Figure \ref{fig_pareto}. Here, the red, blue and green squares correspond to the homotopy runs with $\delta=0$, i.e., they represent the Pareto fronts for the different bi-objective optimization problems (see also the two-dimensional projections in Figure \ref{fig_pareto}). The points obtained by $\delta=0.1$ are represented by the larger circles (again red, blue and green correspond to $\tilde{\mathcal H}^{1,2}$, $\tilde{\mathcal H}^{2,3}$ and $\tilde{\mathcal H}^{3,1}$, respectively), followed by those for $\delta=0.2$ and $\delta=0.3$ in the center. Black symbols denote the transitions between two different homotopies for a given value of $\delta$, i.e., they are solutions of both $\tilde{\mathcal H}^{i,i+1}(\Omega, 1, \delta)=0$ and $\tilde{\mathcal H}^{i+1,i+2}(\Omega, 0, \delta)=0$ (where indices are modulo 3). The black squares represent the minimizers of the three single-objective optimization problems. The shapes corresponding to these points as well as the solution to $\tilde{\mathcal H}^{1,2}(\Omega, 2/3, 0.2)=0$ are shown in Figure \ref{fig_pareto}. We again used a mesh consisting of 650 triangles and 384 vertices. The total computation time for obtaining the depicted Pareto front was about 35 minutes which underlines the efficiency of the method compared to conventional methods.

\newcommand{\paretoPath}{output_bu_20240409/}

%%%%%%%%%%%%%%%%%%%%%%%%%%%%%%%%%%%%%%%%%%%%%%%%%%%%%%%%%%%%
\section{Conclusion and Outlook} \label{sec_conclusion}
In this work, we applied the concept of homotopy methods to (free-form) shape optimization problems without and with PDE constraints. We followed a predictor-corrector approach where the corrector stage is realized by an unregularized shape-Newton method. We discussed how to obtain arbitrary order predictors in shape optimization and investigated their effects numerically. Moreover, we discussed different strategies to automatically adapt the step sizes in the homotopy based on properties of the physical problem. We saw that homotopy methods allow for solving shape optimization problems to high accuracy in a Newton-type manner even with initial guesses that are far away from the optimal solution.
The presented approach could be further extended in the following directions:
\begin{itemize}
 \item  The mentioned challenges concerning mesh quality (cf. Remark \ref{rem_meshqual}) could be circumvented when coupling the approach with a level set method and remeshing routines \cite{AllaireDapogny2021}, thus also allowing for the solution of topology optimization problems in a Newton-type manner.
 \item For physical problems defined on more complicated geometries, the proposed simple PDE-constrained problem \eqref{eq_defG_pde_concrete} may not be feasible. Thus, different types of homotopies such as the Keller model could be investigated.
 \item Often, mathematical models of physical problems involve parameters for which experience shows that gradually updating them is beneficial for the solution process. Examples include incremental loading in nonlinear elasticity or the degree of penalization of intermediate material densities in density-based topology optimization.
 When using a homotopy method, these parameters could either be defined as functions of the homotopy parameter, or could be treated as independent, additional homotopy parameters, amounting to a multi-dimensional homotopy method; see also \cite{Hillermeier2001JOTA,Hillermeier2001book} for the case of multi-objective optimization. In combination with adaptive step size choices, updates of hard-to-tune parameters may be computed automatically in order to find the best path $\{(t(s), x(t(s))): s \in [0,1], t(0)=0 \in \R^p, t(1)=1 \in \R^p \}$ from the initial design at $t(0)=0\in \R^p$ to the solution of the problem at $t(1)=1 \in \R^p$, $p>1$.
 \item The mesh size of the discretized domain could be considered as an additional continuation parameter when coupled with a posteriori error estimators, meaning to start with a coarse mesh and gradually refining the mesh as one approaches the problem of interest; see e.g. \cite{Bandara2016} or \cite{Evgrafov2014} for related ideas.
 \item The concept of continuation methods allows the user to influence the path that is taken from the initial design to a (usually only locally) optimal design. This aspect can be exploited to systematically explore the design space and even detect different paths leading to different local solutions by means of deflation techniques \cite{Papadopoulos2021}.
 \item In the context of multi-objective shape optimiization, the concept of homotopy methods seems to be a promising tool in practice as it allows the designer to interactively steer the design into a desired direction, e.g. in order to fulfill certain manufacturability requirements.
 \item Finally, throughout this paper, we assumed the intermediate solutions to be parametrized by the homotopy parameter and that no turning points appear. While this assumption was fulfilled in the considered examples, it may not hold in general. Therefore, the presented approach could be extended to the case where the solution curve is parametrized by arc length \cite{allgower2012numerical} allowing for a wider class of problems to be treated.
\end{itemize}

\section*{Acknowledgments}
The work of A.C. and P.G. is supported by the FWF funded project P32911 as well as
the joint DFG/FWF Collaborative Research Centre CREATOR (CRC – TRR361/F90) at TU Darmstadt,
TU Graz, RICAM and JKU Linz. The work of B.E. is funded by a Humboldt Postdoctoral
Fellowship and has been supported by
the Cluster of Excellence PhoenixD (EXC 2122, Project ID 390833453). The support is greatly acknowledged.

\bibliographystyle{abbrv}
\bibliography{lit,litPG}

\begin{thebibliography}{10}

\bibitem{AllaireDapogny2021}
{\sc G.~Allaire, C.~Dapogny, and F.~Jouve}, {\em Chapter 1 - shape and topology
  optimization}, in Geometric Partial Differential Equations - Part II,
  A.~Bonito and R.~H. Nochetto, eds., vol.~22 of Handbook of Numerical
  Analysis, Elsevier, 2021, p.~1–132.

\bibitem{allgower2012numerical}
{\sc E.~L. Allgower and K.~Georg}, {\em Numerical continuation methods: an
  introduction}, vol.~13, Springer Science \& Business Media, 2012.

\bibitem{Bandara2016}
{\sc K.~Bandara, T.~Rüberg, and F.~Cirak}, {\em Shape optimisation with
  multiresolution subdivision surfaces and immersed finite elements}, Computer
  Methods in Applied Mechanics and Engineering, 300 (2016), p.~510–539.

\bibitem{BoltenDoganayGottschalkKlamroth2021}
{\sc M.~Bolten, O.~T. Doganay, H.~Gottschalk, and K.~Klamroth}, {\em Tracing
  locally {P}areto-optimal points by numerical integration}, SIAM Journal on
  Control and Optimization, 59 (2021), pp.~3302--3328.

\bibitem{shapeHomotopyZenodo2024}
{\sc A.~Cesarano, B.~Endtmayer, and P.~Gangl}, {\em Supplementary material to
  "{H}omotopy methods for higher order shape optimization: A globalized
  shape-{N}ewton method and {P}areto-front tracing}, 2024.
\newblock https://zenodo.org/doi/10.5281/zenodo.11108761.

\bibitem{CesaranoGanglSCEE}
{\sc A.~Cesarano and P.~Gangl}, {\em Tracing pareto-optimal points for
  multi-objective shape optimization applied to electric machines}, 2024.
\newblock https://arxiv.org/abs/2404.12205.

\bibitem{DeckelnickHerbertHinze2022}
{\sc K.~Deckelnick, P.~J. Herbert, and M.~Hinze}, {\em A novel {$W^{1,\infty}$}
  approach to shape optimisation with {L}ipschitz domains}, ESAIM: COCV, 28
  (2022), p.~2.

\bibitem{b_DEZO_2011a}
{\sc M.~C. Delfour and J.~P. Zol{\'{e}}sio}, {\em Shapes and geometries},
  Society for Industrial and Applied Mathematics, 2011.

\bibitem{Delfour1991}
{\sc M.~C. Delfour and J.~P. Zolésio}, {\em Anatomy of the shape {H}essian},
  Annali di Matematica Pura ed Applicata, 159 (1991), p.~315–339.

\bibitem{Deuflhard2011}
{\sc P.~Deuflhard}, {\em Newton Methods for Nonlinear Problems}, vol.~35 of
  Springer Series in Computational Mathematics, Springer Berlin Heidelberg,
  2011.

\bibitem{Dunlavy2005}
{\sc D.~M. Dunlavy and D.~P. O'Leary}, {\em Homotopy optimization methods for
  global optimization.}

\bibitem{EtlingHerzog2018}
{\sc T.~Etling and R.~Herzog}, {\em Optimum experimental design by shape
  optimization of specimens in linear elasticity}, {SIAM} Journal of Applied
  Mathematics, 78 (2018), pp.~1553--1576.

\bibitem{EtlingHerzogLoayzaWachsmuth2020}
{\sc T.~Etling, R.~Herzog, E.~Loayza, and G.~Wachsmuth}, {\em First and second
  order shape optimization based on restricted mesh deformations}, SIAM Journal
  on Scientific Computing, 42 (2020), pp.~A1200--A1225.

\bibitem{Evgrafov2014}
{\sc A.~Evgrafov}, {\em State space {N}ewton’s method for topology
  optimization}, Computer Methods in Applied Mechanics and Engineering, 278
  (2014), p.~272–290.

\bibitem{Feppon2019Sep}
{\sc F.~Feppon, G.~Allaire, F.~Bordeu, J.~Cortial, and C.~Dapogny}, {\em {Shape
  optimization of a coupled thermal fluid-structure problem in a level set mesh
  evolution framework}}, SeMA, 76 (2019), pp.~413--458.

\bibitem{FEPPON2021113638}
{\sc F.~Feppon, G.~Allaire, C.~Dapogny, and P.~Jolivet}, {\em Body-fitted
  topology optimization of 2d and 3d fluid-to-fluid heat exchangers}, Computer
  Methods in Applied Mechanics and Engineering, 376 (2021), p.~113638.

\bibitem{GLLMS2015}
{\sc P.~Gangl, U.~Langer, A.~Laurain, H.~Meftahi, and K.~Sturm}, {\em Shape
  optimization of an electric motor subject to nonlinear magnetostatics}, SIAM
  Journal on Scientific Computing, 37 (2015), pp.~B1002--B1025.

\bibitem{GanglEtAlSAMO2020}
{\sc P.~Gangl, K.~Sturm, M.~Neunteufel, and J.~Sch\"{o}berl}, {\em Fully and
  semi-automated shape differentiation in {NGS}olve}, Structural and
  Multidisciplinary Optimization, 63 (2020), p.~1579–1607.

\bibitem{Hillermeier2001JOTA}
{\sc C.~Hillermeier}, {\em Generalized homotopy approach to multiobjective
  optimization}, Journal of Optimization Theory and Applications, 110 (2001),
  p.~557–583.

\bibitem{Hillermeier2001book}
{\sc C.~Hillermeier}, {\em Nonlinear Multiobjective Optimization: A Generalized
  Homotopy Approach}, Birkh\"{a}user Basel, 2001.

\bibitem{Houta2023pre}
{\sc Z.~Houta, F.~Messine, and T.~Huguet}, {\em {Topology Optimization for
  Magnetic Circuits with Adjoint Method in 3D}}.
\newblock working paper or preprint, May 2023.

\bibitem{a_IGSTWE_2018a}
{\sc J.~A. Iglesias, K.~Sturm, and F.~Wechsung}, {\em Two-dimensional shape
  optimization with nearly conformal transformations}, {SIAM} Journal on
  Scientific Computing, 40 (2018), pp.~A3807--A3830.

\bibitem{KovtunenkoOhtsuka2022}
{\sc V.~A. Kovtunenko and K.~Ohtsuka}, {\em Inverse problem of shape
  identification from boundary measurement for {S}tokes equations: Shape
  differentiability of {L}agrangian}, Journal of Inverse and Ill-posed
  Problems, 30 (2022), pp.~461--474.

\bibitem{Malinen2010}
{\sc I.~Malinen and J.~Tanskanen}, {\em Homotopy parameter bounding in
  increasing the robustness of homotopy continuation methods in multiplicity
  studies}, Computers \& Chemical Engineering, 34 (2010), p.~1761–1774.

\bibitem{MartinSchuetze2018}
{\sc A.~Martin and O.~Sch\"utze}, {\em Pareto tracer: a predictor–corrector
  method for multi-objective optimization problems}, Engineering Optimization,
  50 (2018), p.~516–536.

\bibitem{Mueller2021}
{\sc P.~M. M\"{u}ller, N.~K\"{u}hl, M.~Siebenborn, K.~Deckelnick, M.~Hinze, and
  T.~Rung}, {\em A novel p-harmonic descent approach applied to fluid dynamic
  shape optimization}, Structural and Multidisciplinary Optimization, 64
  (2021), p.~3489–3503.

\bibitem{NovruziPierre2002}
{\sc A.~Novruzi and M.~Pierre}, {\em Structure of shape derivatives}, Journal
  of Evolution Equations, 2 (2002), p.~365–382.

\bibitem{Onyshkevych2021}
{\sc S.~Onyshkevych and M.~Siebenborn}, {\em Mesh quality preserving shape
  optimization using nonlinear extension operators}, Journal of Optimization
  Theory and Applications, 189 (2021), p.~291–316.

\bibitem{Papadopoulos2021}
{\sc I.~P.~A. Papadopoulos, P.~E. Farrell, and T.~M. Surowiec}, {\em Computing
  multiple solutions of topology optimization problems}, SIAM Journal on
  Scientific Computing, 43 (2021), pp.~A1555--A1582.

\bibitem{Pinzon2022}
{\sc J.~Pinzon and M.~Siebenborn}, {\em Fluid dynamic shape optimization using
  self-adapting nonlinear extension operators with multigrid preconditioners},
  Optimization and Engineering, 24 (2022), p.~1089–1113.

\bibitem{a_SC_2018a}
{\sc S.~Schmidt}, {\em Weak and strong form shape {H}essians and their
  automatic generation}, {SIAM} Journal on Scientific Computing, 40 (2018),
  pp.~C210--C233.

\bibitem{schmidt2008pareto}
{\sc S.~Schmidt and V.~Schulz}, {\em Pareto-curve continuation in
  multi-objective optimization}, Pacific Journal of Optimization, 4 (2008),
  pp.~243--258.

\bibitem{SchmidtSchulz2023}
{\sc S.~Schmidt and V.~H. Schulz}, {\em A linear view on shape optimization},
  SIAM Journal on Control and Optimization, 61 (2023), pp.~2358--2378.

\bibitem{SchmidthWadbroBerggren2016}
{\sc S.~Schmidt, E.~Wadbro, and M.~Berggren}, {\em Large-scale
  three-dimensional acoustic horn optimization}, SIAM Journal on Scientific
  Computing, 38 (2016), pp.~B917--B940.

\bibitem{Schoeberl2014}
{\sc J.~Sch\"oberl}, {\em C++11 implementation of finite elements in
  {NGSolve}}, Tech. Rep.~30, Institute for Analysis and Scientific Computing,
  Vienna University of Technology, 2014.

\bibitem{SchulzSiebenbornWelker2016}
{\sc V.~H. Schulz, M.~Siebenborn, and K.~Welker}, {\em Efficient {PDE}
  constrained shape optimization based on {S}teklov–{P}oincaré-type
  metrics}, SIAM Journal on Optimization, 26 (2016), p.~2800–2819.

\bibitem{Simon1989}
{\sc J.~Simon}, {\em Second variations for domain optimization problems}, in
  Control and estimation of distributed parameter systems (Vorau, 1988),
  vol.~91 of Internat. Ser. Numer. Math., Birkh\"auser, Basel, 1989,
  p.~361–378.

\bibitem{SZ}
{\sc J.~Soko{\l}owski and J.-P. Zol{\'e}sio}, {\em Introduction to shape
  optimization}, vol.~16 of Springer Series in Computational Mathematics,
  Springer-Verlag, Berlin, 1992.
\newblock Shape sensitivity analysis.

\bibitem{a_ST_2015a}
{\sc K.~Sturm}, {\em Minimax {L}agrangian approach to the differentiability of
  nonlinear {PDE} constrained shape functions without saddle point assumption},
  {SIAM} Journal on Control and Optimization, 53 (2015), pp.~2017--2039.

\bibitem{sturm2018convergence}
{\sc K.~Sturm}, {\em Convergence of {N}ewton's method in shape optimisation via
  approximate normal functions}, 2018.

\bibitem{suli2003introduction}
{\sc E.~S{\"u}li and D.~F. Mayers}, {\em An introduction to numerical
  analysis}, Cambridge university press, 2003.

\bibitem{WatsonHaftka1989}
{\sc L.~T. Watson and R.~T. Haftka}, {\em Modern homotopy methods in
  optimization}, Computer Methods in Applied Mechanics and Engineering, 74
  (1989), pp.~289--305.

\bibitem{Zulehner1988}
{\sc W.~Zulehner}, {\em A simple homotopy method for determining all isolated
  solutions to polynomial systems}, Mathematics of Computation, 50 (1988),
  pp.~167--177.

\end{thebibliography}

%
% \newpage
% {\magenta
% TODOlist:
% \textbf{TODOs}
% \begin{itemize}
%     \item[$\checkmark$] P: Abstract
%     \item[$\checkmark$] make connection to EtlingHerzogLoayzaWachsmuth2020 more clear
%     \item[$\checkmark$] Bernhard step secant predictors shorten
%     \item[$\checkmark$] Upload code, make bib entry
%     \item[$\checkmark$] Conclusion
%     \item[$\checkmark$] remark mesh quality
%     \item[$\checkmark$] remark cost higher order predictors
%     \item[$\checkmark$] non-connected parts of pareto
%     \item[$\checkmark$] re-think name ``agile'', ``double agile''; also in numerics section (even in titles of plots)!
%      \item[$\checkmark$] cite Zulehner;
%     \item[$\checkmark$] explain why $d \mathcal J(\Omega)(n)$ bzw explain vector $d \mathcal J(\Omega)(V_x^{normal})$, for all boundary points $x$; also in Algorithm 3.1
%     \item[$\checkmark$] allow larger step sizes for gradient method
%     \item[$\checkmark$] Maybe remove tang reg and BFGS? If so, adjust text everywhere (remove discussion on tang and bfgs; will get shorter). If not, add reference for SciPy in footnote and finish ``may even outperform on finer mesh''...?
%     \item[$\checkmark$] re-run shape grad for ex9
%     \end{itemize}
% }

%%%%%%%%%%%%%%%%%%%%%%%%%%%%%%%%%%%%%%%%%%%%%%%%%%%%%%%%%%%

\end{document}